\documentclass[10pt,letterpaper]{article}

\usepackage[utf8]{inputenc}
\usepackage[T1]{fontenc}

\usepackage{graphicx,epsf}
\usepackage{graphbox}
\usepackage{bbm,colo	r}
\usepackage{amsmath,amssymb,amsfonts}
\usepackage{multirow}
\usepackage{bm}
\usepackage{algorithmic}
\usepackage[ruled,vlined,linesnumbered]{algorithm2e}
\usepackage{hyperref}
\usepackage{amsthm}
\usepackage{dsfont}
\usepackage{subfig}
\usepackage{nccmath}
\usepackage{color}

\usepackage[left=2.8cm, right=2.8cm, top=2.5cm, bottom=2.5cm]{geometry}

\usepackage{cite}

\newcommand{\R}{\mathbb{R}}
\newcommand{\C}{\mathbb{C}}
\newcommand{\M}{\mathbb{M}}
\newcommand{\HH}{\mathcal{H}}

\newtheorem{theorem}{Theorem}[section]

\graphicspath{
{./Figures_Final/}
}

\title{Optimization Methods for Joint Eigendecomposition}

\author{Erik Troedsson\footnote{Centre for Mathematical Sciences, Lund University, Lund, Sweden, (erik.troedsson@math.lu.se)}, Marcus Carlsson\footnote{Centre for Mathematical Sciences, Lund University, Lund, Sweden, (marcus.carlsson@math.lu.se)}, Herwig Wendt\footnote{IRIT laboratory, Univ. Toulouse, CNRS, Toulouse, France, (herwig.wendt@irit.fr)}}

\begin{document}

\maketitle

\begin{abstract}
Joint diagonalization, the process of finding a shared set of approximate eigenvectors for a collection of matrices, arises in diverse applications such as multidimensional harmonic analysis or quantum information theory. This task is typically framed as an optimization problem: minimizing a non-convex function that quantifies off-diagonal matrix elements across possible bases. In this work, we introduce a suite of efficient algorithms designed to locate local minimizers of this functional. Our methods leverage the Hessian’s structure to bypass direct computation of second-order derivatives, evaluating it as either an operator or bilinear form -- a strategy that remains computationally feasible even for large-scale applications. Additionally, we demonstrate that this Hessian-based information enables precise estimation of parameters, such as step-size, in first-order optimization techniques like Gradient Descent and Conjugate Gradient, and the design of second-order methods such as (Quasi-)Newton. The resulting algorithms for joint diagonalization outperform existing techniques, and we provide comprehensive numerical evidence of their superior performance.

\begin{center}\bf Keywords\end{center} matrix diagonalization, simultaneous diagonalization, joint eigen-decomposition, gradient descent, conjugate gradient, Newton method
\end{abstract}


\section{Introduction}

\noindent{\bf Context and motivation.\quad}
Let $\mathcal{A} = \{A_1,\ldots,A_K\}$ be a collection of $n\times n$ real or complex matrices. Joint eigen-decomposition, often referred to as joint diagonalization (JD), is the process of finding an invertible matrix $U$ that approximately diagonalizes all matrices in $\mathcal{A}$.
This technique is crucial in a variety of applications across multiple disciplines, including blind source separation \cite{jutten1991blind,cardoso1996equivariant,pham2001joint,comon2010handbook}, independent component analysis \cite{cardoso1993blind,comon1994independent,ziehe2004fast}, common principal component problem \cite{flury1984common}, canonical polyadic decomposition \cite{roemer2013semi,luciani2014canonical}, decoupling of VAR models \cite{xu2024testing}, multidimensional harmonic analysis \cite{haardt2008higher,andersson2015method,sahnoun2017multidimensional,andersson2018esprit}, geometric multimodal learning \cite{behmanesh2021geometric}, or quantum information theory \cite{burel2022joint}. In these contexts, JD enables the extraction of meaningful underlying structures from complex datasets by identifying approximate eigenvector sets shared across matrix collections.
For a matrix $U$ to jointly diagonalize each $A_k\in\mathcal{A}$, the transformed matrices $U^{-1}A_kU$ should ideally be diagonal. In that case, $U$ can be explicitly found by doing an eigen-decomposition of either matrix $A_k$. However, in practical scenarios, perfect JD is rarely achievable due to noise or other data inconsistencies. Thus, the focus shifts to finding an approximate solution that minimizes the off-diagonal components of each transformed matrix.

 \noindent{\bf Problem formulation.\quad}
 The JD problem is commonly formulated as an optimization problem.
 Specifically, the objective is to minimize a functional $ f_\mathcal{A}(U)$ over all  invertible matrices $U$, where $ f_\mathcal{A}(U)$ measures the collective size of the off-diagonal elements of each $U^{-1}A_kU$,
\begin{equation}\label{functional}
    f_\mathcal{A}(U) = \frac{1}{2}\sum_{k=1}^K\sum_{i\neq j} \left|(U^{-1}A_kU)_{ij}\right|^2.
\end{equation}
This function, which is non-convex and challenging to minimize globally, represents the squared sum of the off-diagonal elements in the transformed matrices $U^{-1}A_kU$.
The goal is to identify local minimizers of $f_\mathcal{A}(U)$, which correspond to approximate joint diagonalizers for $\mathcal{A}$.

\noindent{\bf Related work and limitations.\quad}
Numerous approaches have been proposed for solving the JD problem. The problem has often been considered in combination with additional assumptions such as $U^{-1}=U^*$, which naturally appears if the matrices $A_k$ are Hermitian, see, e.g., \cite{pham2001blind,yeredor2002non,ziehe2004fast,tichavsky2008fast,souloumiac2009nonorthogonal,ablin2018faster,azzouz2023generalized} and references therein. 
Here we focus on the general, unconstrained JD problem, which has received lesser attention.
Some methods apply LU decompositions \cite{luciani2015joint},  while others leverage manifold optimization techniques \cite{mesloub2018efficient,iferroudjene2009new,gong2012complex}, or randomized algorithms \cite{he2024randomized}.
Other more recent methods simplify $f_{\mathcal{A}}$ using Taylor series approximations \cite{andre2020joint,cao2022joint}, to name a few. Despite these advances, existing algorithms often fail to achieve true minimization of $f_{\mathcal{A}}$, especially in high-dimensional settings or when handling complex data structures.

A key aspect of our approach involves  a closed form formula for the gradient and Hessian of the objective function $f_{\mathcal{A}}$.
Although the gradient has been discussed in earlier works \cite{hori1999joint}, it has, to our knowledge, not been used in a practical and operational algorithm.
Additionally, our method leverages on the Hessian of $f_{\mathcal{A}}$, expressed as a bilinear form rather than as a full matrix. This is a core innovation that enables efficient computation and avoids scaling issues for large matries.

\noindent{\bf Contributions.\quad}
This paper presents a set of novel algorithms for JD that address these limitations and advance the state of the art. The primary contributions of this work are:
\begin{enumerate}
\item {\it Efficient Hessian Computation}: We develop a technique to evaluate the Hessian of $f_{\mathcal{A}}$ as either a bilinear form or a linear operator. This approach avoids the direct computation and storage of second-order partial derivatives, enabling scalable application to large matrices.
\item {\it Theoretical Insights into Objective Function Properties}: We present new theoretical findings on the structure of the objective function $f_{\mathcal{A}}$.
Particularly, we show that it diverges at non-invertible points, allowing us to set safe step-size bounds that prevent convergence issues.We also show that a multiplicative change of basis improves the descent path, enhancing convergence speed.
\item {\it Optimization Schemes Using Second Order Information}: By leveraging the Hessian's structure, we introduce an effective step-size selection method for classical first-order optimization techniques, including Gradient Descent and Conjugate Gradient. We further adapt second-order optimization schemes, such as Quasi-Newton methods,  to improve convergence.
\item {\it New Gradient-Based Algorithms}: We propose innovative algorithms that integrate these insights, outperforming existing JD methods in terms of both computational efficiency and solution quality. These algorithms demonstrate robust performance, as shown in our numerical experiments.
\end{enumerate}
Our algorithms are designed to operate seamlessly on both complex and real-valued data. 
For concreteness, we focus here on the unconstrained case, but they are also easily adaptable to the Hermitian case.

\noindent{\bf Outline.\quad}
The remainder of this paper is organized as follows. Section \ref{sec:grad} provides a formal problem formulation, along with explicit expressions for the gradient and Hessian of $f_{\mathcal{A}}$. Section \ref{sec:theory} discusses the theoretical properties of $f_{\mathcal{A}}$, including its non-convexity and the implications of non-invertible points. In Section \ref{sec:algo}, we introduce our proposed algorithms, detailing the integration of the gradient, Hessian, and optimized step-size selection in our descent methods.
Section \ref{sec:results} presents a comprehensive set of numerical experiments that benchmark our algorithms against state-of-the-art JD methods. Finally, Section \ref{sec:3Dharmonic} explores an application to a 3D harmonic retrieval problem, and Section \ref{sec:conclusions} concludes with insights and potential avenues for future work.

\section{{Problem formulation, Gradient and Hessian}}
\label{sec:grad}

\subsection{Problem formulation and notation}
Let $\M_n(\R)$ and $\M_n(\C)$ denote the set of real and complex matrices of dimension $n\times n$, where we simply write $\M_n$ if it is of no importance whether they are real or complex.
Let $\mathcal{A} = \{A_1,\ldots,A_K\}$ be a collection of $K$ matrices in $\M_n$, let $J$ denote the matrix with zeros on the diagonal and ones elsewhere and let $\diamond$ denote Hadamard-multiplication of matrices. With this new notation the functional \eqref{functional} can be written
\begin{equation}\label{functional2}
    f_\mathcal{A}(U) = \frac{1}{2}\sum_{k=1}^K
     \left|J \diamond(U^{-1}A_kU)\right|^2_F,
\end{equation}
where $F$ is the Frobenius norm. The joint diagonalizer of $\mathcal{A}$ is defined as the matrix (or possibly matrices) $\hat U$ that solves the problem
\begin{equation}\label{problem}
\hat U = \arg\min_{U}    f_\mathcal{A}(U).
\end{equation}
In the following, we will define and study descent algorithms for solving \eqref{problem} that make use of new, explicit expressions of the gradient and Hessian of the objective function $ f_\mathcal{A}(U)$. Note that the joint eigen-vectors can be complex-valued even if $\mathcal{A}$ is real. To see this just consider the trivial case $K=1$ and the matrix
$A_1=\left(                                                                                                                                                                  \begin{array}{rr}                                                                                                                                                  0 & 1 \\
-1 & 0 \\
\end{array}                                                                                                                                                                \right)$
whose eigenvectors are $[1,\pm i]^t$ (where $t$ denotes matrix transpose). Thus, we are obliged to consider $f_{\mathcal{A}}$ as a real-valued functional on $\M_n(\C)$.

\subsection{Gradients and Hessians for real-valued functions on the complex linear space $\M_n(\C)$.}

Before getting to the formulas for the gradient and bilinear Hessian, we need to clarify how these are to be realized only relying on ``natural'' operations in $\M_n(\C)$, that is, avoiding vectorization and partial derivatives. Indeed, the classical approach would be to consider $\M_n(\C)$ as a real linear space of dimension $2n^2$ and then, upon introducing a concrete basis, compute $\nabla f_{\mathcal{A}}$ as a vector in $\R^{2n^2}$ and the Hessian as a matrix of dimensions $2n^2\times 2n^2$, by computing roughly $(2n^2)^2/2$ partial derivatives. However, this is not computationally feasible, and we will proceed to introduce basis free versions of the above objects. First of all, we note that $\M_n(\C)$ can be considered as an inner product space over $\R$ by introducing a new real-valued inner product (a.k.a.~scalar product) $\langle U,V \rangle$ simply by defining it as $\mathsf{Re} \langle U,V\rangle $ where the latter refers to the complex Frobenius inner product,
$$
\langle U,V\rangle =\sum_{i,j}U_{i,j}\overline{V}_{i,j}.
$$ 
We also remind the reader that a symmetric bilinear form $\HH$ on $\M_n(\C)$ is a function on $\M_n(\C)\times \M_n(\C)$ which is linear in each argument separately and also satisfies $\HH(U,V)=\HH(V,U)$ for all $U,V\in\M_n(\C)$.

Armed with this, we now rely on linear algebra to show that, for each fixed $U_0\in\M_n(\C)$, there is a unique element $\nabla f|_{U_0}\in \M_n(\C)$ and a unique symmetric bilinear form $\mathcal{H}|_{U_0}$ such that $f(U_0+Z)$ can be written as
\begin{multline}\label{cf1}
f(U_0+Z)=f(U_0)+\mathsf{Re}\langle \nabla f|_{U_0}(Z),Z\rangle\\+\frac{1}{2}\HH|_{U_0}(Z,\!Z)+\mathcal{O}(\|Z\|^3)\end{multline}
where $\mathcal{O}$ denotes the so called ``big ordo'', we refer to \cite{troedsson2024joint} for the details. Moreover, given any bilinear form $\HH$, there is always a unique (real) linear operator $H:\M_n(\C)\rightarrow \M_n(\C)$ such that \begin{equation}\label{cf2}\HH(U,V)=\mathsf{Re}\langle H(U),V\rangle.\end{equation}
Clearly this operator is symmetric if and only if the bilinear form is symmetric.

\subsection{Gradient and Hessian for the functional \eqref{functional2}}
\label{sec:GradHess}
The previous section indicates that we can obtain coordinate free versions of the gradient and of the classical Hessian matrix by simply finding rules to compute objects that satisfy the above criteria. For the concrete functional $f_{\mathcal{A}}$ we provide formulas in this section but refer to the parallel paper \cite{troedsson2024joint} for proofs.
In Section \ref{sec:algo}, we will utilize them to construct gradient descent, conjugate gradient and Quasi-Newton algorithms, and to provide step-size selection rules for them. Let $U^*$ be the conjugate transpose of a given matrix $U$, and let $U^{-*}$ be shorthand notation for $(U^{-1})^*$.
\begin{theorem}\label{t1}
    The gradient and bilinear Hessian of \eqref{functional2} are
\begin{align}
        &\nabla f_\mathcal{A}|_U = \sum_{k=1}^K  U^{-*} \left[D_k^*,D_k\diamond J\right], \label{gradient}\\
        &\HH_{\mathcal{A}}(Z,W)=\sum_{k=1}^K\mathsf{Re}\langle J \diamond [D_k, U^{-1}Z],[D_k, U^{-1}W] \rangle \label{hessianbilinear} \\&+ \mathsf{Re}\langle J\diamond D_k, [U^{-1}Z, U^{-1}WD_k]+[U^{-1}W, U^{-1}ZD_k]\rangle,\nonumber
            \end{align}%
where $D_k\triangleq U^{-1}A_kU$.
\end{theorem}
The proofs are found in the publication \cite{troedsson2024joint}, see Theorems 5.1 and 6.1. The main idea is to use Taylor expansions, e.g.,
$$(U+Z)^{-1}=\sum_{k=0}^2 (U^{-1}Z)^k U^{-1}+\mathcal{O}(\|Z\|^3)$$
and then simply multiply terms together, collecting first order terms separately, second order terms separately, and to rely on the uniqueness of \eqref{cf1}.

To arrive at the linear operator incarnation of the Hessian involves more work, basically consisting in reshuffling terms in the expression for $\HH_{\mathcal{A}}(Z,W)$ to isolate $W$ on one side of the inner product. The resulting expression is:
\begin{theorem}\label{t11}
   The Hessian operator of \eqref{functional2} is given by
\begin{align}&H_\mathcal{A}|_U(Z)=\sum_{k=1}^K {U^{-*}}\Big(\big[D_k^*,J\diamond [D_k,U^{-1}Z]\big]+\label{hessian}\\ &\big[(U^{-1}Z)^*,J\diamond D_k\big]D_k^* +\left[J\diamond D_k,\left(U^{-1}ZD_k\right)^*\right]\Big).\nonumber\end{align}
\end{theorem}
While this is efficiently implemented as a forward operator, it is, in contrast to matrix incarnations, not something which is straightforward to invert. Despite this, we shall see in Section \ref{sec:QN} that it is nonetheless possible to design efficient and competitive Quasi-Newton methods based on this expression.

\smallskip

\noindent{\bf Complexity analysis.\quad}
Assuming that standard matrix multiplication and the Gauss-Jordan algorithm for computing the matrix inverse are used, the complexity of computing the gradient \eqref{gradient}, bilinear Hessian \eqref{hessianbilinear} of Hessian operator \eqref{hessian} is $\mathcal{O}(Kn^3)$, respectively. This is up to constants the same complexity as an evaluation of the objective function \eqref{functional2}. The memory requirement is $\mathcal{O}(Kn^2)$.
Note that the use of the Hessian in matrix form would, in contrast, require the storage of $n^4$ elements, and its computation, inversion, and application $\mathcal{O}(Kn^4)$, $\mathcal{O}(n^6)$ and $\mathcal{O}(n^4)$ operations, respectively.

\subsection{Descents algorithms with multiplicative change of basis}

A basic descent algorithm for solving \eqref{problem} starts with an initial invertible matrix $U_0$ and then iteratively uses a descent direction $S_{\mathcal{A}}|_{U_0}$ computed at this point, for instance the gradient $S_\mathcal{A}|_{U_0}=-\nabla f_\mathcal{A}|_{U_0}$ or some approximation thereof, as search direction to compute the next point as
\begin{equation}\label{standardupdate}
    U_{m+1} = U_m + \lambda_m S_\mathcal{A}|_{U_m}
\end{equation}
where $\lambda_m >0$ is a step size parameter. In this setting the matrix $U_m$ is updated additively.
An alternative update can be defined by observing that the objective function in \eqref{problem} has the property
\begin{equation}\label{ut6}
f_\mathcal{A}(U_m V)=f_{U_m^{-1} \mathcal{A} U_m}(V).
\end{equation}
Therefore, one can instead consider computing the descent direction at each iteration for the modified objective function $f_{U_m^{-1} \mathcal{A} U_m}$ evaluated at the identity $I$, and update via
\begin{equation}
\label{multupdate}
U_{m+1}=U_m(I+\lambda_mS_{U_m^{-1} \mathcal{A} U_m}|_I).
\end{equation}
Updates for descent algorithms with such a multiplicative change of basis for \eqref{problem} have first been studied in \cite{cardoso1996equivariant} and subsequently been used in, e.g, \cite{andre2020joint}, 
and empirically observed to yield better results \cite{ErikEUSIPCO2024}.
In Section \ref{sec:multiplicative}, we will study the properties of such updates and compare them with those of standard additive descent updates.
Numerically, our results also indicate that this alternative significantly improves convergence speed of our proposed descent algorithms (see Sec.~\ref{sec:result:update}). Moreover, the gradient \eqref{gradient} and the Hessian (\ref{hessianbilinear}-\ref{hessian}) are less expensive to compute at the identity matrix because no matrix inverse is involved (since $U=I$ and $I^{-1}=I$).

To provide some intuition why the use of a multiplicative change of basis is beneficial, the local geometry of the graph of $f_{\mathcal{A}}$ is illustrated in the top row of Fig.~\ref{fig:objfun}
whereas the local geometry of the corresponding graph $f_{U^{-1}\mathcal{A}U}$ is shown in the bottom row. Note that this is only a two dimensional slice of a $2n^2$-dimensional space, where one direction is given by the normalized gradient, and the other is just a randomly chosen orthonormal direction.
As is plain to see, the graph for $f_{U^{-1}\mathcal{A}U}$ looks much less chaotic than that for $f_{\mathcal{A}}$, and the singularities (i.e.~points where the function approaches $\infty$) are further away from the current iterate. In the next section, we provide some theoretical remarks concerning $f_{\mathcal{A}}$ which aim to partially clarify these observations, and theoretical bounds that allow us to ensure that we avoid stepping near the singularities in the descent algorithms that we propose in Section \ref{sec:algo}.

\section{Theoretical properties of $f_{\mathcal{A}}$}
\label{sec:theory}

\subsection{Non-convexity and singularities of $f_{\mathcal{A}}$}
\label{sec:nonconvex}

The functional $f_{\mathcal{A}}$ is highly non-convex, note that even without noise or data inconsistencies it has more than one basin of attraction because of its permutation indeterminacy. 
Moreover, it is not even defined for non-invertible matrices. Despite this fact, nothing a priori prevents the iterates of an algorithm to converge towards such a point, i.e., it is not clear whether such points could form local minimizers on the graph of $f_{\mathcal{A}}$. However, the key result of the companion paper \cite{troedsson2024joint} states that this cannot happen, with the exception of some degenerate scenarios which happen with probability zero. The result reads as follows:

\begin{theorem}
  \label{kern}
   Let $U$ denote any invertible matrix. We have that
    \begin{equation}\label{blowup}\lim_{U\to U_0}f_{\mathcal{A}}(U) = \infty\end{equation}
    for all fixed nonzero non-invertible matrices $U_0$, if and only if the matrices in $\mathcal{A}$ have no common nontrivial invariant subspace.
\end{theorem}
In other words, when the condition is fulfilled, the singular points of $f_{\mathcal{A}}$ are precisely the non-invertible matrices that one would like to avoid. An illustration of this fact is provided Fig.~\ref{fig:objfun:sing}; note in particular that every time $\lambda$ hits a value where the corresponding matrix is non-invertible, the graph goes to $\infty$, as stipulated by the theorem.
To phrase it differently, applying descent algorithms to $f_{\mathcal{A}}$ as a means to find a complete set of joint eigenvectors, is a well posed strategy.
The proof of this theorem is quite intricate, the details are provided in \cite[Sec.~2]{troedsson2024joint}. In the same paper \cite[Sec.~3]{troedsson2024joint} it is shown that, whenever the matrices $\mathcal{A}$ are corrupted by noise drawn from an absolutely continuous probability measure, it holds with probability one that the matrices in $\mathcal{A}$ have no common nontrivial invariant subspace.
Therefore, {in most applications}, we can be sure that local minima of \eqref{functional2} are invertible. 

With the above facts in mind, an appropriate choice of step size to approximately minimize $f_\mathcal{A}$ over the line $U+\lambda S$, given a descent direction $S$, becomes critical. The following result provides a theoretical upper bound for reasonable step-lengths. We remind the reader that the spectral radius $\rho(S)$ is the size of the modulus of the largest eigenvector of a given matrix $S$.

\begin{figure}[tb]
\centering
\includegraphics[width=0.33\linewidth,trim=0 22 20 10, clip]{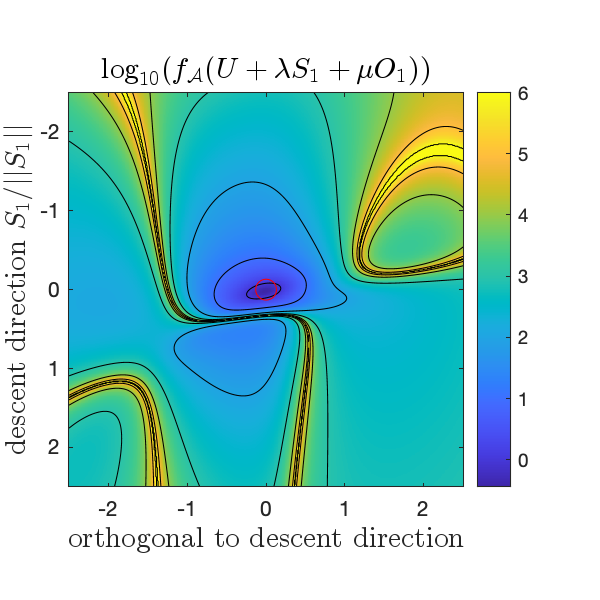}%
\includegraphics[width=0.33\linewidth,trim=0 22 20 10, clip]{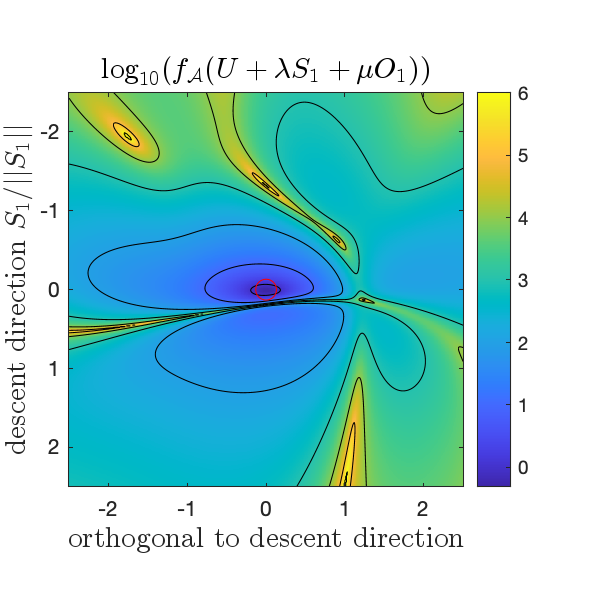}%
\\
\includegraphics[width=0.33\linewidth,trim=0 22 20 10, clip]{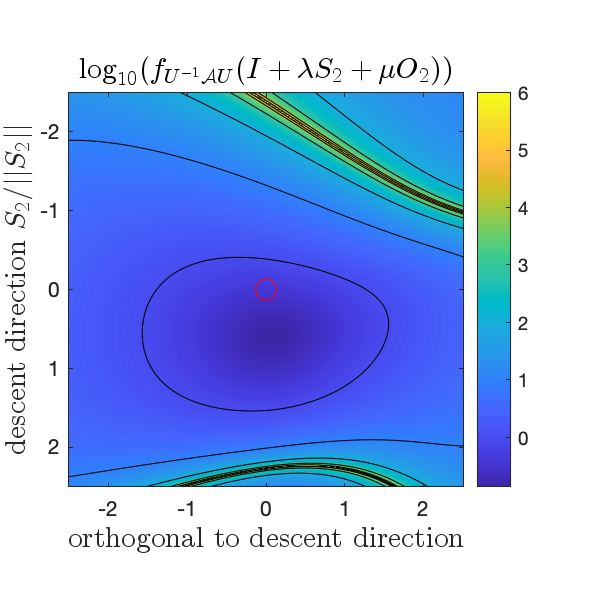}%
\includegraphics[width=0.33\linewidth,trim=0 22 20 10, clip]{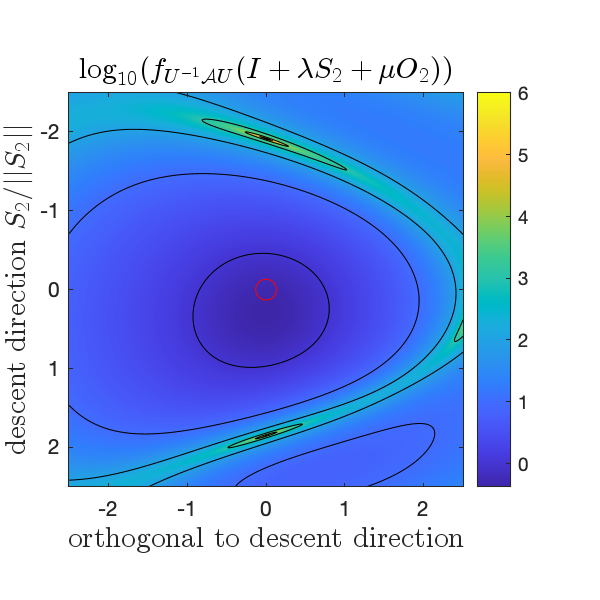}%
\caption{\label{fig:objfun}
Top: The objective function $f_\mathcal{A}$ \eqref{functional2} for a matrix tuple $\mathcal{A}$, constructed as described in Section \ref{sec:MC} with SNR=$30$dB, with dimensions $n=10,K=5$ (left) and $n=10,K=10$ (right) plotted along the line $U+x S_1/\|S_1\| + y O_1$, where $S_1$ is the gradient \eqref{equ:grad:S1} and $O_1$ is a random orthogonal direction with norm $1$. The center plots show the functional $f_{U^{-1}\mathcal{A}U}$ plotted along the line $I+x S_2/\|S_2\|+ y O_2$, where $S_2$ is the gradient \eqref{equ:grad:S2} and $O_2$ is a random orthogonal direction with norm $1$.
In both cases $U=\mbox{eig}(\sum_{k=1}^K A_k)$,  which is a natural choice for initialization for the JD problem (indicated by the red circle in the graphs).}
\end{figure}

\begin{theorem}\label{t2}
If $U$ is invertible, the objective function  $f_\mathcal{A}(U+\lambda S)$ has no singularity on the interval \begin{equation}\label{h6}
                                                                                                                0\leq \lambda<\frac{1}{\rho(U^{-1}S)}.
                                                                                                              \end{equation}  Similarly, the objective function  $f_{U^{-1}\mathcal{A}U}(I+\lambda S)$ has no singularity on the interval $0\leq \lambda<1/\rho(S)$.
\end{theorem}

\indent{\it Proof.\quad} Considering the expression \eqref{functional2} for the objective function $f_\mathcal{A}(U)$, it is clear that $\lambda\mapsto f_\mathcal{A}(U+\lambda S)$ may have poles at the singularities of $(U+\lambda S)^{-1}$. Since $(U+\lambda S)^{-1}=U^{-1}(I+\lambda U^{-1}S)$ one sees that these poles are precisely minus the reciprocal of the non-zero eigenvalues of $U^{-1}S$. Therefore, the first potential singularity of $f_\mathcal{A}(U+\lambda S)$ is at distance $1/\rho(U^{-1}S)$ from the origin. The second part follows upon substitution of $U$ by $I$. $\square$

Note that the estimate based on the spectral radius is sharp, i.e., there exists a $\lambda$ on the boundary where the matrix is not invertible. In practice, it is numerically more efficient to use $1/\|U^{-1}S\|$ as an upper bound in \eqref{h6}, since the norm is faster to compute and we clearly have that $\rho(M)\leq \|M\|$ for any matrix $M$;
see Sec.~\ref{sec:stepsize} for details on our proposed step size selection strategy. 
In Fig.~\ref{fig:objfun:sing}, the bounds $1/\rho(U^{-1}S_1)$ and $1/\rho(S_2)$ are plotted and indeed coincide with the position of the first non-invertible point encountered along the descent directions, respectively.

\subsection{Properties of multiplicative change of basis}
\label{sec:multiplicative}

When using the multiplicative change of basis, the objective function $f_{U_m^{-1}\mathcal{A}U_m}$ is different at each iteration, but its critical points are related with those of the original objective function $f_{\mathcal{A}}$ through a linear change of variables, see \eqref{ut6}. Here we will study how this change of variables affects the descent algorithms, with the goal to understand the empirical observation that updates with multiplicative change of basis work better.
Again, we refer to Fig.~\ref{fig:objfun} which seems to indicate that the graph of $f_{U_m^{-1}\mathcal{A}U_m}$ is more stable and suitable for optimization, than the original version $f_{\mathcal{A}}$. See also Fig.~\ref{fig:basis}, which indicates a massive improvement on algorithmic performance.

Let us assume that we are at a specific matrix $U_m=U$, which is close to a local minimizer of $f_\mathcal{A}$, and compare a standard gradient step with a gradient step with multiplicative change of basis.
First, making use of \eqref{gradient}, the direction in the multiplicative update \eqref{multupdate} is given by 
\begin{equation}\label{equ:grad:S2}
S_2\triangleq-\nabla f_{U^{-1}\mathcal{A}U}|_{I}=-\sum_{k=1}^K  \left[D_k^*,(D_k\diamond J)\right],
\end{equation}
whereas the direction in the ``standard update'' \eqref{standardupdate} is given by
\begin{equation}\label{equ:grad:S1}
S_1\triangleq-\nabla f_{\mathcal{A}}|_U= U^{-*}S_2. 
\end{equation}
Note here that if $U$ is close to diagonalizing $\mathcal{A}$, $D_k=U^{-1}A_kU$ is likely a ``well-behaved'' matrix which is close to diagonal. In particular, $D_k\diamond J$ will be close to zero, and hence it is also to be expected that $S_2$ is a small matrix, in some sense.
In the standard update, the objective function we want to minimize, as a function of the step-length $\lambda$, is
\begin{align*}
     & f_{\mathcal{A}}(U+\lambda S_1)=f_{\mathcal{A}}(U+\lambda U^{-*}S_2)=\\
&\nonumber \sum_{k=1}^K\left\|J\diamond \Big((U+\lambda U^{-*}S_2)^{-1}A_k(U+\lambda U^{-*}S_2)\Big)\right\|_F=\\
     &\nonumber \sum_{k=1}^K\left\|J\diamond \Big((I+\lambda U^{-1}U^{-*}S_2)^{-1}D_k(I+\lambda U^{-1}U^{-*}S_2)\Big)\right\|_F.
   \end{align*}
In contrast, when using gradient descent with multiplicative change of basis as in \eqref{multupdate}, the corresponding functional is
\begin{align*}
  f_{U^{-1}\!\mathcal{A}U}(I\!+\!\lambda S_2)\!=\! 
  \sum_{k=1}^K\left\|J\!\diamond\! \Big((I+\lambda S_2)^{-1}D_k(I+\lambda S_2)\Big)\right\|_F.
   \end{align*}
By comparing the two expressions, we see that the only difference is the factor $U^{-1}U^{-*}$ in front of $S_2$ in the expression for the standard update. 
The matrix $U^{-1}U^{-*}$ is Hermitian and thus has a singular value decomposition $V^*\Sigma V$ where $V$ is unitary and $\Sigma$ is diagonal. Since the Frobenius norm is invariant under multiplication by unitary matrices, we can think of the operation $S_2\mapsto U^{-1}U^{-*}S_2$ as taking the matrix $S_2$ and first appling an isometry (multiplication by $V$), followed by the multiplication of $\Sigma$ and finally another isometry. But  in realistic applications, it is quite likely that $\Sigma$ has both small and large singular values. In other words, $\Sigma$ may have a large condition-number, and thus the multiplication with $\Sigma$  will elongate certain directions in the matrix space, and compress others. 

\begin{figure}[tb]
\centering
\scriptsize $n=10,K=5$\hskip24mm $n=10,K=10$\vskip-0mm
\includegraphics[width=0.45\linewidth]{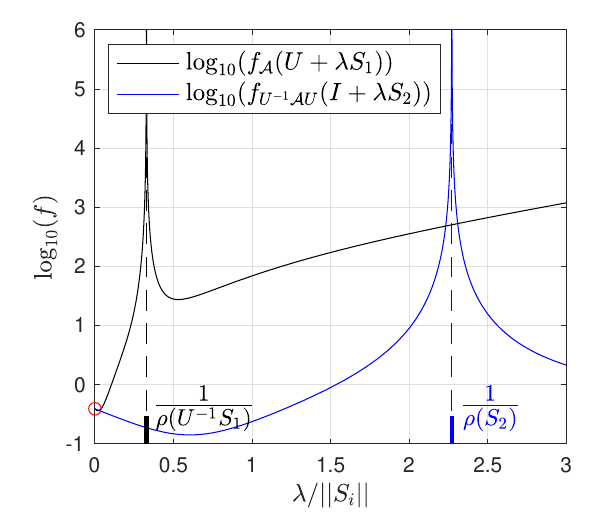}%
\includegraphics[width=0.45\linewidth]{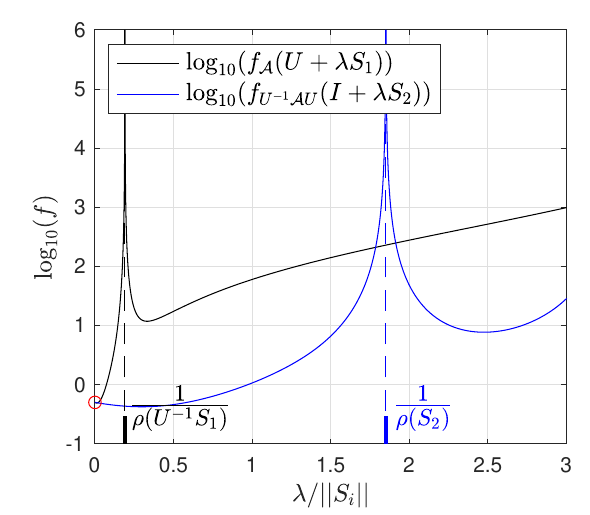}
\caption{\label{fig:objfun:sing}
Slices of the functionals $f_\mathcal{A}$ and $f_{U^{-1}\mathcal{A}U}$, as simulated and depicted in Fig.~\ref{fig:objfun}, along the descent directions $S_1$ and $S_2$ respectively. The dashed vertical lines indicate values of $\lambda$ for which $U+\lambda S_1$ (black) and $I+\lambda S_2$ (blue) are non-invertible, respectively; the bold solid lines at the bottom indicate the corresponding prediction given by Theorem \ref{t2}.
}
\end{figure}

\noindent{\bf Numerical illustration.\quad}While the outlined argument is not entirely rigorous, it correlates very well with what is demonstrated in Fig.~\ref{fig:objfun}, where the graph of $f_{\mathcal{A}}$ displays very different behavior in the two (orthogonal and normalized) directions. In this example, $\mbox{cond}(U^{-1}U^{-*})$ equals $86.2$ (left subplot, $173.1$ for the right), which can explain why the upper plots in Figure \ref{fig:objfun} are more chaotic than the lower ones. 
Similarly, Fig.~\ref{fig:objfun:sing} shows that the spectral radius of $U^{-1}S_1=U^{-1}U^{-*}S_2$ is much larger than that of $S_2$ - and the possible function value decrease with one step along $S_1$ is much smaller than for $S_2$: more precisely,
we have that $1/\rho(U^{-1}S_1)$ equals $0.33$ (left subplot, $0.19$ for the right) but $2.27$ ($1.85$) for $1/\rho(S_2)$, and the function value decrease equals $-0.026$ ($-0.012$) for $S_1$, while it equals $-0.223$ ($-0.078$) for $S_2$, respectively.

\section{Descent algorithms with multiplicative basis change}
\label{sec:algo}

Armed with the results and insights from the previous sections, we can now detail our proposed algorithms for solving the JD problem.

\subsection{Gradient descent algorithm}
\label{sec:gd}

Given Theorem \ref{t1}, the gradient descent is straightforward to implement.
Setting $A_{k,m}=U_m^{-1}A_k U_m$ and writing $\mathcal{A}_{m}$ for $\{A_{1,m},\ldots,A_{K,m}\}$,
the use of the gradient as a search direction, combined with the multiplicative change of basis, leads to the following scheme:
\begin{enumerate}
\item[] {\it Initialization} ($m=0$): set $\mathcal{A}_{0}=\mathcal{A}$ and $U_0=U_{init}$.
    \item Given $\mathcal{A}_{m}$, compute the gradient $\nabla f_{\mathcal{A}_{m}}|_I$.
    \item Compute the step size $\lambda_m$.
    \item Compute $\mathcal{A}_{m+1}$ by setting $$A_{k,m+1} = (I-\lambda_{m} \nabla f_{\mathcal{A}_{m}}|_I)^{-1} A_{k,m} (I-\lambda_{m} \nabla f_{\mathcal{A}_{m}}|_I).$$
    \item Update the eigenvectors to $$U_{m+1} = U_{m}(I-\lambda_{m}\nabla f_{\mathcal{A}_{m}}|_I).$$
    \item If not converged, set $m=m+1$ and go back to (i).
\end{enumerate}
A natural choice for initializing the matrix $U$ is given by
\begin{equation}\label{equ:Uinit}
U_{init}=\mbox{eig}\left(\sum_{k=1}^K A_k\right),
\end{equation}
see \cite{Yeredor2005} for a closely related alternative; assuming absence of prior knowledge, one could also consider, e.g., $U_{init}=I$.
An effective rule for determining the step size $\lambda_m$ in step (ii)  will be given in Section \ref{sec:stepsize} (cf., \eqref{equ:finalstep}). Note that if only the eigenvalues, but not the eigenvectors $U$, are of interest, step (iv) can be omitted.

\subsection{Conjugate gradient descent algorithm}
\label{sec:cgd}

Standard gradient descent is a simple and robust method but unfortunately it can have slow convergence, and it is well known that the conjugate gradient method often is faster.
However, since the matrices $\mathcal{A}_{m}$ are changed at each iteration, we are in each step considering a different functional $f_{\mathcal{A}_m}$, and therefore the usual scheme has to be modified.
Recall that the idea of the conjugate gradient method in the standard setting (with fixed matrices $\mathcal{A}_{m}=\mathcal{A}$) is to replace the gradient directions $-\nabla f_{\mathcal{A}}(U_m)$ with search directions $S_m$ and ensure that these directions are orthogonal with respect to the inner product induced by the Hessian operator \eqref{hessian}
\begin{equation}\label{conjugate}
    \langle S_{m}, H_\mathcal{A}|_{U_{m}}(S_{m-1})\rangle =0,
\end{equation}
for all $m$ \cite{daniel1967conjugate}. In the multiplicative update setting, we end up with a curious situation where, after the matrix tuple $\mathcal{A}_{m-1}$ is updated to $\mathcal{A}_{m}$, we are effectively at the point $U_{m}=I$ and are now seeking to minimize the new functional $f_{\mathcal{A}_{m}}$, where \begin{equation}\label{Aupdate}
                                                                                                      \mathcal{A}_{m}=(I+\lambda_{m-1} S_{m-1})^{-1}\mathcal{A}_{m-1} (I+\lambda_{m-1} S_{m-1}).
                                                                                                    \end{equation}
The above transformation of $\mathcal{A}_{m-1}$ corresponds to a change of basis on the underlying matrix space we are optimizing over, where objects $V$ in step $m-1$ get transformed by $(I+\lambda_{m-1} S_{m-1})^{-1}V$. Thus, the old direction $S_{m-1}$ will, with respect to the new tuple $\mathcal{A}_{m}$, correspond to
\begin{equation}\label{oldsearchinnewbases}
    \tilde{S}_{m-1} \triangleq (I+\lambda_{m-1} S_{m-1})^{-1}S_{m-1}.
\end{equation}
We now take the new search direction $S_{m}$ to be of the form
\begin{equation}
\label{form}
S_{m} = -\nabla f_{\mathcal{A}_{m}}|_I + \beta_m \tilde{S}_{m-1}.
\end{equation}
Inserting the above into the conjugate criterion \eqref{conjugate}, we obtain
$$0\!=\!-\langle  \nabla f_{\mathcal{A}_{m}}|_I,H_{\mathcal{A}_{m}}|_I(\tilde{S}_{m-1})\rangle + \beta_m \langle \tilde{S}_{m-1},H_{\mathcal{A}_{m}}|_I(\tilde{S}_{m-1})\rangle,$$
and thus we have the following formula for $\beta_m$,
\begin{equation}\label{beta}
    \beta_m = \frac{\langle  \nabla f_{\mathcal{A}_{m}}|_I,H_{\mathcal{A}_{m}}|_I(\tilde{S}_{m-1})\rangle}{\langle \tilde{S}_{m-1},H_{\mathcal{A}_{m}}|_I(\tilde{S}_{m-1})\rangle}.
\end{equation}
Without \eqref{oldsearchinnewbases} owed to the multiplicative change of basis, this is known as Daniel's rule for nonlinear conjugate gradient \cite{daniel1967conjugate}, which fundamentally differs from the more commonly used Fletcher-Reeves, Dai-Yuan, etc. formulas for $\beta$ \cite{dai1999nonlinear, fletcher1964function, hestenes1952methods, polak1969note}, in that it makes explicit use of the Hessian for second order information from $f_\mathcal{A}$. The main critique of Daniel's rule is that, at first sight, it seems to require the evaluation of the full Hessian matrix \cite{william2006survey}. Indeed, one can in principle extract a matrix representation of the Hessian by identifying $\mathbb{M}(n,\mathbb{C})$ with the real space $\mathbb{R}^{2n^2}$, introducing an orthonormal basis and computing the second order partial derivatives.
 However, the use of the matrix increases the problem size quadratically, so this approach does not scale well and quickly becomes numerically intractable.
This highlights the value of Theorem \ref{t1}, which allows us to efficiently compute the Hessian for the directions of interest in \eqref{beta}, whilst avoiding costly storage and application of matrix incarnations of the Hessian.

The proposed conjugate gradient descent method hence consists of the following steps:
\begin{enumerate}
\item[] {\it Initialization} ($m=0$): Set $\mathcal{A}_{0}=\mathcal{A}$, $U_0=U_{init}$ and $S_{0}=-\nabla f_{\mathcal{A}}|_I$ and perform a gradient descent step.
    \item Given $S_{m-1}$ and $\mathcal{A}_{m}$, compute $\tilde{S}_{m-1}$ via \eqref{oldsearchinnewbases}.
    \item Compute the gradient $\nabla f_{\mathcal{A}_{m}}|_I$ and the new search direction $S_{m}$ according to \eqref{form} and \eqref{beta}.
    \item Compute the  step-size $\lambda_{m}$.
    \item Compute $\mathcal{A}_{m+1}$ via \eqref{Aupdate} (using $m+1$ in place of $m$).
    \item Update the eigenvectors to $U_{m+1} = U_{m}(I+\lambda_{m}S_{m})$.
    \item If not converged, set $m=m+1$ and go back to (i).
\end{enumerate}
Again, we will present an effective step size rule below and step (v) is unnecessary if only eigenvalues are sought.

The conjugate gradient method is proven to converge in a number of convex but also non-convex instances, see e.g.~\cite{dai1999nonlinear,gilbert1992global,william2006survey} and the references therein. In particular, our numerical results show that the proposed algorithm works well for the non-convex functional \eqref{functional2} (see Section \ref{sec:results}). However, for non-convex problems, the Hessian can take negative values at a given step, which is a strong indication that the conjugate gradient direction may be completely off. On the same note, it is known that convergence is improved if we discard negative $\beta$-values, see in particular \cite{gilbert1992global}. In summary, if $\beta_m$ is negative or if the denominator in \eqref{beta} is negative, we simply set $\beta_m=0$ so that $S_m$ becomes a gradient descent step.

\begin{figure}[htb]
\centering
\includegraphics[width=0.45\linewidth]{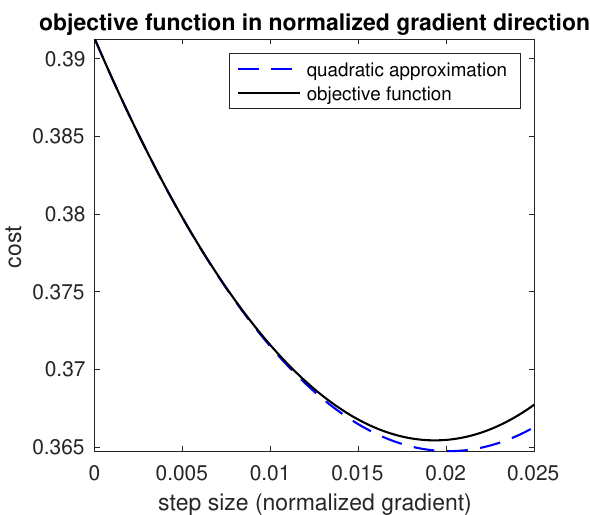}%
\includegraphics[width=0.45\linewidth]{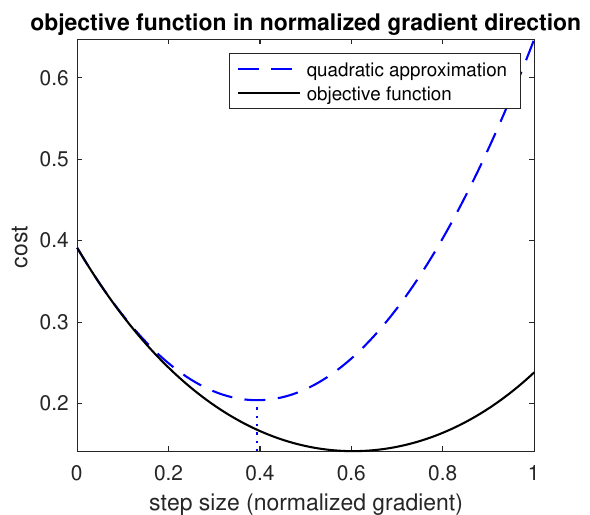}%
\caption{\label{fig:step1}{\bf Step size selection rule.} Objective function (black solid) and approximation \eqref{equ:approxStep} (blue dashed) in the gradient direction for the additive (left) and multiplicative (right) update corresponding to the left column of Fig.~\ref{fig:objfun:sing}.}
\end{figure}

\subsection{Newton method and Quasi-Newton method}
\label{sec:QN}
The Newton method is a classical descent algorithm that makes explicit use of second order information, i.e., it relies on knowing the Hessian, in order to improve the convergence rate. However, like in the previous section, we demonstrate here that it is not necessary to have access to the Hessian {matrix}, the operator version given by Theorem \ref{t11} is sufficient. Specifically, in Newton's method, given a point $U$, the search direction $S$ is given by the solution to
\begin{equation}
	H_\mathcal{A}|_U(S) = -\nabla f_{\mathcal{A}}(U). \label{newt}
\end{equation}
A standard method to solve \eqref{newt} would be to vectorize both sides and apply a suitable solver for linear equations. This would however require us to form the matrix version of the Hessian. To avoid this large computation, we instead propose to solve \eqref{newt} by minimizing the quadratic objective function
\begin{equation}\label{q}
	q_{\mathcal{A},U}(S) = \frac{1}{2}\mathsf{Re}\langle S,H_\mathcal{A}|_U(S)\rangle + \mathsf{Re}\langle S,\nabla f_{\mathcal{A}}(U)\rangle,
\end{equation}
for which we will only need the operator version \eqref{hessian} of the Hessian. Here it is crucial to take the real part of the inner product, since we view $\mathbb{M}(n,\mathbb{C})$ as a real inner product space, see \cite{troedsson2024joint} for further details.
Note that $\nabla q_{\mathcal{A},U}|_S = H_\mathcal{A}|_U(S)+\nabla f_{\mathcal{A}}(U)$,
so any stationary point of \eqref{q} solves \eqref{newt}.
To minimize \eqref{q}, the standard conjugate gradient algorithm can be used in a straight-forward manner, as detailed in the pseudocode given below.

Our proposed approach comes with several advantages. Firstly, we only need to apply $H_\mathcal{A}|_U$ to matrices (which is possible since we have the operator form), instead of computing the large $2n^2\times 2n^2$ representation. Secondly, we may terminate this conjugate gradient algorithm early to get an approximate solution of \eqref{newt}, which yields a trade-off between finding an exact solution and reduced computation time. The resulting ``Quasi-Newton'' method thus consists of an outer loop of the descent algorithm for $f_\mathcal{A}$,
and an inner loop to minimize $q_{\mathcal{A},U}$, where $U=I$ when using a multiplicative change of basis, and is summarized below.
\begin{enumerate}
	\item[] {\it Initialization} ($m=0$): Set $\mathcal{A}_0 = \mathcal{A}$ and $U_0=U_{init}$.
	\item Given $\mathcal{A}_{m}$ and $U_{m}$, compute $S_{m}$ by minimizing \eqref{q}:
	\begin{enumerate}
	\item[] {\it Initialization} ($l=0$):\\ $s_{0}=\nabla f_{\mathcal{A}|I}$, $r_0=s_0-H_{\mathcal{A}|I}(s_0)$, $p_0=r_0$
	\item $\alpha_l=\frac{\mathsf{Re}\langle r_l,r_l\rangle}{\mathsf{Re}\langle r_l,H_{\mathcal{A}|I}(r_l)\rangle}$
	\item $s_{l+1}=s_l+\alpha_lp_l$, $r_{l+1}=r_l-\alpha_lH_{\mathcal{A}|I}(p_l)$
	\item $\beta_l=\frac{\mathsf{Re}\langle r_{l+1},r_{l+1}\rangle}{\mathsf{Re}\langle r_{l},r_{l}\rangle}$
	\item $p_{l+1}=r_{l+1}+\beta_lp_{l}$
	\item if converged set $S_m=s_{l+1}$ and go to (ii), \\
	else set $l=l+1$ and go back to (a).
	\end{enumerate}
	\item Compute the  step-size $\lambda_{m}$.
	\item Compute $\mathcal{A}_{m+1}$ via \eqref{Aupdate} (using $m+1$ in place of $m$).
    \item Update the eigenvectors to $U_{m+1} = U_{m}(I+\lambda_{m}S_{m})$.
    \item If not converged, set $m=m+1$ and go back to (i).
\end{enumerate}
In the numerical simulations reported in Section \ref{sec:results}, the stopping criterion used in step (e) of the inner loop is that either $||H_\mathcal{A}|_U(S) +\nabla f_{\mathcal{A}}(I)||_F^2$ has fallen below 10\% of its initial value, or a maximum number of $l=100$ iterations have been reached.
Again, we will present an effective step size rule below and step (iv) is unnecessary if only eigenvalues are sought. 
Finally, when the Hessian $ H_\mathcal{A}|_U$ is not positive semi-definite, the minimization of \eqref{q} may fail since then $q_{\mathcal{A},U}$ does not have a minimum. In this case, we use the negative gradient as the search direction instead.

\begin{figure}[tb]
\centering
\includegraphics[width=0.4\linewidth]{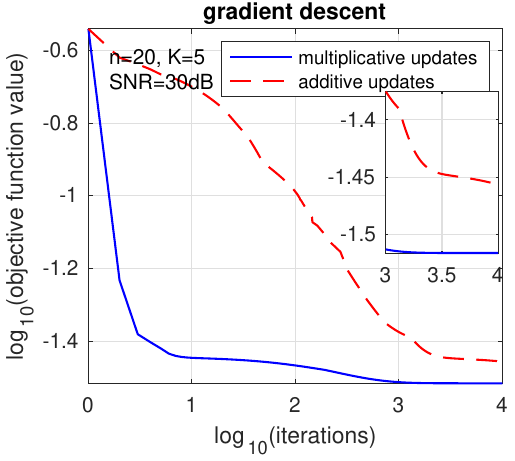}%
\includegraphics[width=0.4\linewidth]{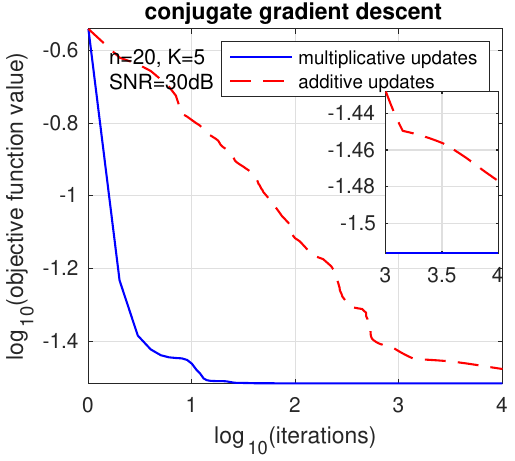}
\caption{\label{fig:basis}{\bf Additive vs. multiplicative updates.\quad}Median of the objective function values versus iterations with additive \eqref{standardupdate} and multiplicative \eqref{multupdate} updates  for gradient descent and conjugate gradient descent.
 }
\end{figure}
\begin{figure}[htb]
\centering
\begin{minipage}{0.5\linewidth}
\includegraphics[width=\linewidth,trim=35 0 30 0, clip]{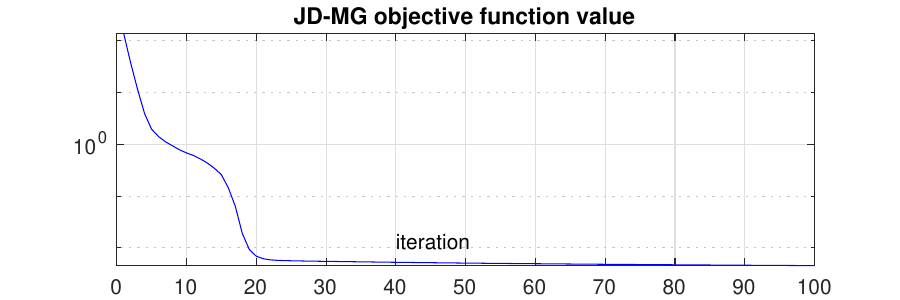}\\
\includegraphics[width=\linewidth,trim=35 0 30 0, clip]{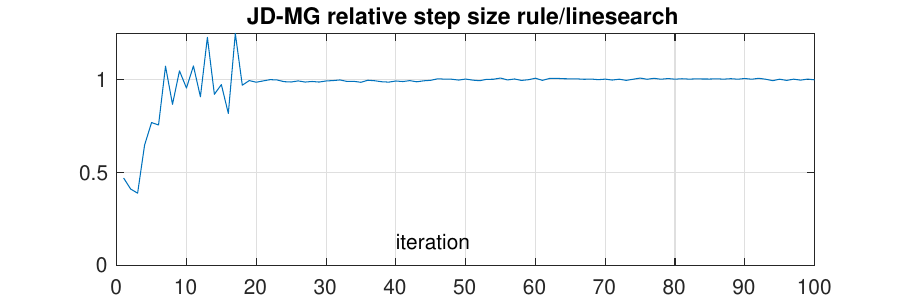}
\end{minipage}%
\begin{minipage}{0.5\linewidth}
\includegraphics[width=0.5\linewidth]{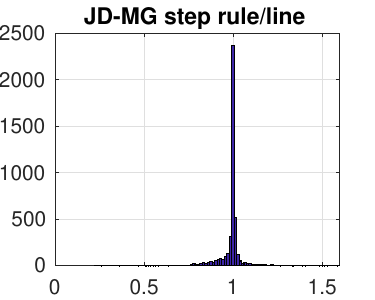}%
\includegraphics[width=0.5\linewidth]{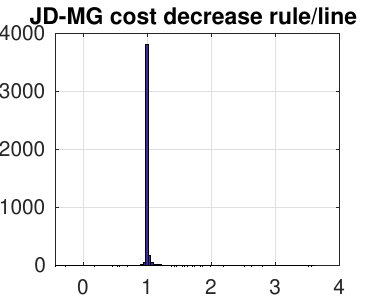}%
\end{minipage}%
\caption{\label{fig:step}{\bf Step size selection performance.} Average (over 100 realizations) objective function decrease with step size \eqref{equ:finalstep} (top left), step size \eqref{equ:finalstep} at each iteration divided by step size obtained with line search (bottom left); histogram of step size \eqref{equ:finalstep} divided by step size obtained with line search (center), histogram of relative cost decrease  (right) when step size \eqref{equ:finalstep} or line search are used (initialization at $U_{init}=I$, SNR=30dB, $n=10,K=5$)}
\end{figure}

\begin{figure*}[htb]
\centering
\includegraphics[width=0.33\linewidth]{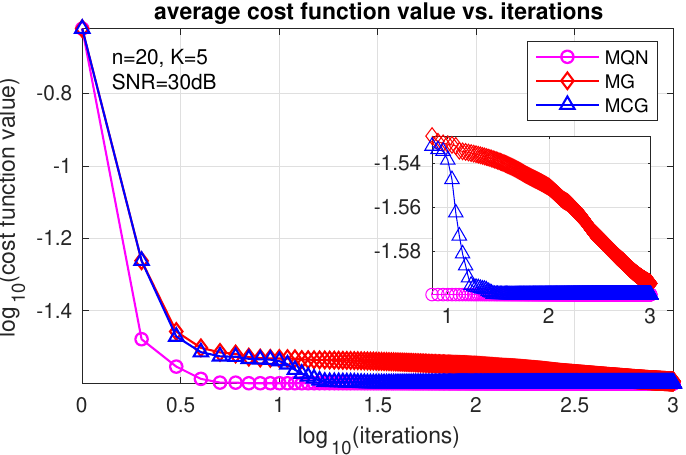}%
\includegraphics[width=0.33\linewidth]{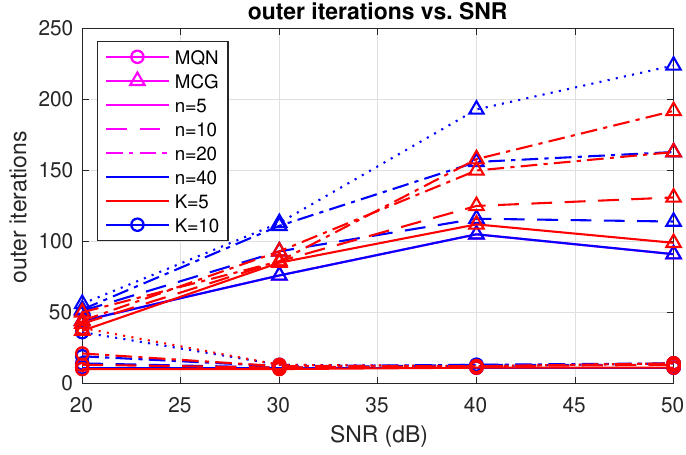}%
\includegraphics[width=0.33\linewidth]{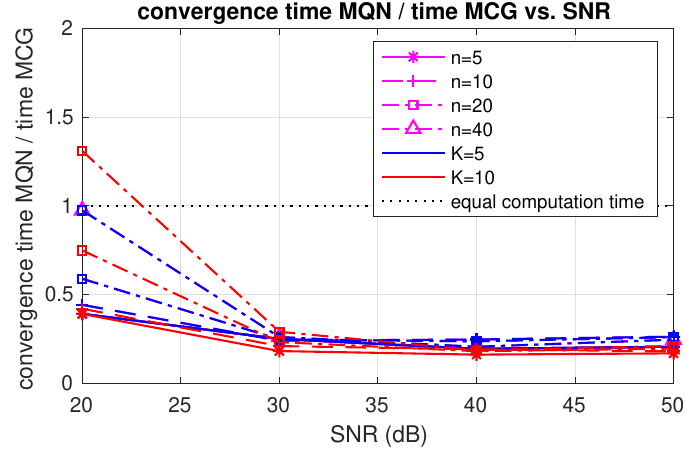}%
\caption{\label{fig:QN}{\bf Performance comparison for proposed descent algorithms.\quad}
Cost function vs. iteration number (left, SNR=$30$dB); Median iteration number to convergence to same objective function value (center); Median of fraction of time to convergence for QN and MCG (right). See Section \ref{sec:perf} for a fuller description.}
\end{figure*}

\subsection{Step size selection rule based on the Hessian}
\label{sec:stepsize}
A standard approach to determining the step size $\lambda_m$ in each iteration of a descent algorithm consists in performing some flavor of line search, i.e., repeated evaluations of the objective function for potential $\lambda_m$'s obtained by some rule (e.g., bisection of a search interval) until a good candidate is found, that is, a candidate that approximately minimize $f_\mathcal{A}$ along the descent direction. This has the disadvantages that it requires to determine a search interval and to perform multiple objective function evaluations at each iteration, which can be costly.
In view of our results in Section \ref{sec:GradHess}, our proposed strategy is to instead make use the bilinear Hessian to find a reasonable step size $\lambda$ given a fixed descent direction $S$ at some point $U$.
Starting with the approximation
\begin{equation}\label{equ:approxStep}
f_{\mathcal{A}}(U+\lambda S)\approx f_\mathcal{A}(U) + \lambda\langle\nabla f_\mathcal{A}|_{U},S\rangle + \frac{1}{2}\lambda^2\langle S,H_{\mathcal{A}}|_{U}(S)\rangle,
\end{equation}
and assuming that $\langle S,H_{\mathcal{A}}|_{U}(S)\rangle>0$, we can minimize \eqref{equ:approxStep} over $\lambda$, which yields the step size rule
\begin{equation}\label{stepsize}
    \lambda_H = -\frac{\langle\nabla f_\mathcal{A}|_{U},S\rangle}{\langle S,H_{\mathcal{A}}|_{U}(S)\rangle}.
\end{equation}
An illustration is given in Fig.~\ref{equ:approxStep}.
There is in general no guarantee for the denominator of \eqref{stepsize} to be positive, and when this is not the case then clearly we can no longer use \eqref{equ:approxStep}. An alternative in this situation is to further study the denominator, which after simplifying can be written as
$$ \sum_{k=1}^K\| J \diamond [D_k, U^{-1}S] \|^2 + 2\langle J\diamond D_k, [U^{-1}S, U^{-1}SD_k]\rangle$$ (recall \eqref{hessianbilinear}).
By discarding the second term inside this sum we obtain a nonnegative Gauss-Newton approximation of the Hessian, which leads to the following step-size rule
\begin{equation}\label{stepsize2}
    \lambda_{GN} = -\frac{\langle\nabla f_\mathcal{A}|_U,S\rangle}{\| J \diamond [D_k, U^{-1}S] \|^2}.
\end{equation}
We propose to use it in cases when the denominator of \eqref{stepsize} is negative.
Moreover, we make use of Theorem \ref{t2} to avoid getting too close to a singularity of $f_{\mathcal{A}}$ and do not allow the step size to be larger than
$$\lambda_{max}=(2\|U^{-1}S\|)^{-1}.$$
Above we use the norm instead of the spectral radius since it is faster to compute, and the factor 2 comes from the fact that we want to stay away from singularities with a good margin, to avoid that the objective function increases.
Eventually, our step size rule is defined by
 \begin{equation}
 \label{equ:finalstep}
 \lambda = \begin{cases}
  \min(\lambda_H,\lambda_{max}) & \textnormal{if }  \langle S,H_{\mathcal{A}}|_{U}(S)\rangle>0\\
  \min(\lambda_{GN},\lambda_{max}) & \textnormal{otherwise.}
 \end{cases}
  \end{equation}
Note that since $f_\mathcal{A}$ is nonconvex, it is possible that the step-size \eqref{equ:finalstep} nevertheless leads to an increase in objective function value.
In this situation, \eqref{equ:finalstep} can nevertheless provide reasonable starting points for performing bisection over the line $U+\lambda S$ to approximately minimize $f_\mathcal{A}$. In our numerical results reported below, we have not made use of this and simply applied rule \eqref{equ:finalstep}.

\section{Numerical Experiments and Results}
\label{sec:results}

\subsection{Monte Carlo simulations}\label{sec:MC}
We study the performance of the proposed suite of algorithms for the JD problem by applying them to a large number of realizations ($N_{MC}=1000$) of collections $\mathcal{A} = \{A_1,\ldots,A_K\}$  of $K$ complex valued matrices of size  $n\times n$, to which we add circular white Gaussian noise for several signal to noise ratios (SNR).
Each matrix $A_k,k=1,\ldots,K$ is constructed as  $A_k=Z\Delta_k Z^{-1}$, where $Z$ is an $n\times n$ matrix whose entries are i.i.d. circular Gaussian random variables, normalized to yield unit norm columns, and $\Delta_k$ are diagonal matrices whose diagonal entries are complex random variables with independent real and imaginary parts that are uniformly distributed on the interval $[-1,1]$.
We denote our proposed  gradient descent (GD), conjugate gradient (CG) and Quasi-Newton algorithms for JD with multiplicative change of basis by JD-MG, JD-MCG and JD-MQN, respectively.
If not mentioned otherwise, all algorithms are initialized at the guess \eqref{equ:Uinit}, and are stopped after $M=1000$ iterations, without other stopping criterion.

\subsection{Multiplicative vs. additive updates}
\label{sec:result:update}

Fig.~\ref{fig:basis} shows the median over realizations of the objective function values versus iterations for GD (left) and CG  (right) when either the standard additive update \eqref{standardupdate} or the multiplicative update  \eqref{multupdate} are used ($n=20,K=5$, SNR=$30$dB, $M=10000$ iterations). The algorithms with multiplicative change of basis clearly converge significantly faster than the ones with standard update: Indeed, the latter require more than one order of magnitude more iterations: while JD-MG and JD-MCG converged after approximately $1000$ and $<100$ iterations, respectively, their counterparts with standard updates have not converged at iteration $10000$. This corroborates our theoretical considerations in Section \ref{sec:multiplicative}.

\subsection{Step size selection performance}
\label{sec:result:step}
To assess the relevance of our step size selection rule \eqref{equ:finalstep}, we apply JD-MG for $100$ realizations (SNR=30dB, $n=10,K=5$, initialization at $U_{init}=I$) and compare, at each step, our step size against that obtained with line search (with 1000 objective function evaluations). The results are plotted in Fig.~\ref{fig:step} and clearly indicate that after a few steps following the initialization (far from the minimum where \eqref{equ:finalstep} yields slightly smaller values), the proposed step size \eqref{equ:finalstep} is essentially identical to the exact minimizer found by line search. Note that in this experiment, line search is 3 orders of magnitude slower than our step size rule. Overall, this leads to conclude that \eqref{equ:finalstep} is an excellent alternative to classical step size selection procedures.

\subsection{Performance comparison for proposed descent algorithms}
\label{sec:perf}
Fig.~\ref{fig:QN} (left) plots median objective function values for our proposed JD-MG, JD-MCG and JD-MQN algorithms ($n=20,K=5$, SNR=$30$dB).  It can be observed that JD-MQN needs fewest iterations (of the order of 10) to reach a minimum, followed by JD-MCG, which needs of the order of 100 iterations to reach the same objective function value, whereas JD-MG still has not converged after 1000 iterations. In the following, we will therefore not consider JD-MG any longer. Fig.~\ref{fig:QN} (center) and Fig.~\ref{fig:QN} (right) provide a more detailed comparison between JD-MCG and JD-MQN; in this experiment, we let JD-MCG converge (i.e., we stop JD-MCG when the relative improvement of the objective function is smaller than $10^{-12}$ times the initial objective function value), stop JD-MQN when it has reached the same objective function value as JD-MCG, and monitor the respective number of iterations and computation times. Fig.~\ref{fig:QN} (center) reports the median number of iterations for various values for $n$, $K$ and SNR. With the exception of very low SNR values, JD-MQN needs significantly less outer iterations than JD-MCG to converge, and dramatically less for large SNR. Fig.~\ref{fig:QN} (right) plots the corresponding relative computation time, showing that JD-MQN can be up to $5$ times faster than JD-MCG for large SNR values. While these results do not exclude that JD-MCG can be a better choice in certain scenarios (low SNR values, bad initial guess $U_{init}$), they nevertheless show that our proposed JD-MQN can be competitive.

\begin{figure}[tb]
\centering
\includegraphics[width=0.45\linewidth]{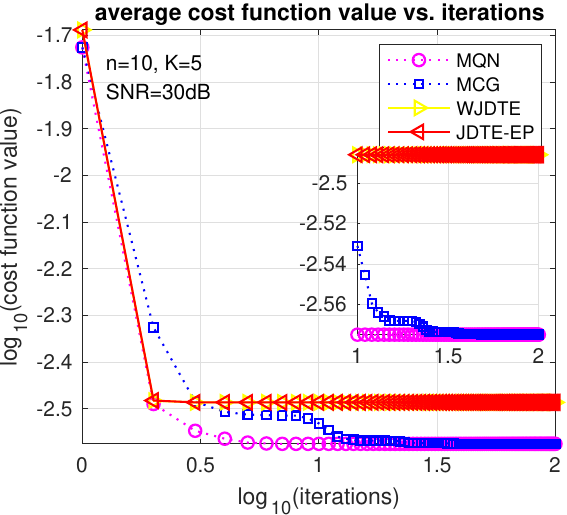}%
\includegraphics[width=0.45\linewidth]{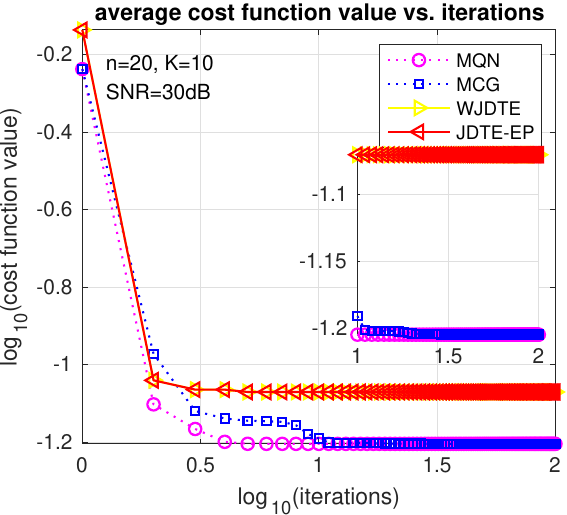}
\caption{\label{fig:perf:iterations}{\bf Performance vs. iterations.\quad}
Median of objective function values versus iterations (SNR=$30$dB;  left column: $n=10,K=5$; right column: $n=20,K=10$).}
\end{figure}

\begin{table*}[tb]
\centering
{\small\setlength{\tabcolsep}{0.1pt}
\begin{tabular}{| r || rrrrrr  || rrrrrr  || rrrrrr |}
\hline
\multicolumn{19} {|c|}{\boldmath$\log_{10}(f_{\mathcal{A}})$}\\
\hline
&\multicolumn{6} {|c||}{$n=10$, $K=5$}
&\multicolumn{6} {|c||}{$n=20$, $K=5$}
&\multicolumn{6} {|c|}{$n=20$, $K=10$}\\
\hline
SNR (dB) 		& $10$ & $20$ & $30$ & $40$ & $50$ & $60$ &
  $10$ & $20$ & $30$ & $40$ & $50$ & $60$ &
    $10$ & $20$ & $30$ & $40$ & $50$ & $60$ 
\\
\hline
eig$_{\sum A}$ &  $1.87$ & $0.29$ & $-1.69$ & $-3.70$ & $-5.70$ & $-7.70$ &
 $2.78$ & $1.29$ & $-0.59$ & $-2.59$ & $-4.59$ & $-6.59$ &
  $3.14$ & $1.73$ & $-0.14$ & $-2.14$ & $-4.14$ & $-6.14$ 
\\
WJDTE 		& $1.23$ & $-0.50$ & $-2.49$ & $-4.48$ & $-6.48$ & $-8.48$ &
 $1.98$ & $0.50$ & $-1.49$ & $-3.49$ & $-5.49$ & $-7.49$ &
  $2.35$ & $0.90$ & $-1.07$ & $-3.07$ & $-5.07$ & $-7.07$ 
\\
JDTE-EP 		& $1.23$ & $-0.50$ & $-2.49$ & $-4.48$ & $-6.48$ & $-8.48$ &
 $1.98$ & $0.50$ & $-1.49$ & $-3.49$ & $-5.49$ & $-7.49$ &
  $2.35$ & $0.90$ & $-1.07$ & $-3.07$ & $-5.07$ & $-7.07$ 
\\
RJD			& $1.48$ & $-0.22$ & $-2.22$ & $-4.22$ & $-6.21$ & $-8.22$ &
 $2.37$ & $0.89$ & $-1.08$ & $-3.08$ & $-5.08$ & $-7.08$ &
  $2.76$ & $1.28$ & $-0.70$ & $-2.69$ & $-4.70$ & $-6.69$ 
\\
JD-MCG 		& $0.99$ & $-0.60$ & $-2.57$ & $-4.57$ & $-6.57$ & $-8.57$ &
 $1.63$ & $0.28$ & $-1.60$ & $-3.60$ & $-5.60$ & $-7.60$ &
  $2.01$ & $0.66$ & $-1.20$ & $-3.20$ & $-5.20$ & $-7.20$ 
\\
JD-MQN 		& $0.99$ & $-0.60$ & $-2.57$ & $-4.57$ & $-6.57$ & $-8.57$ &
 $1.63$ & $0.28$ & $-1.60$ & $-3.60$ & $-5.60$ & $-7.60$ &
  $2.01$ & $0.66$ & $-1.20$ & $-3.20$ & $-5.20$ & $-7.20$ 
\\
\hline
\hline
\multicolumn{19} {|c|}{\boldmath\bf eigenvalue estimation error ($\log_{10}$)}\\
\hline
&\multicolumn{6} {|c||}{$n=10$, $K=5$}
&\multicolumn{6} {|c||}{$n=20$, $K=5$}
&\multicolumn{6} {|c|}{$n=20$, $K=10$}\\
\hline
SNR (dB) 		& $10$ & $20$ & $30$ & $40$ & $50$ & $60$ &
  $10$ & $20$ & $30$ & $40$ & $50$ & $60$ &
    $10$ & $20$ & $30$ & $40$ & $50$ & $60$ 
\\
\hline
eig$_{\sum A}$ & $1.29$ & $-0.67$ & $-3.01$ & $-5.03$ & $-7.03$ & $-9.03$ &
  $1.84$ & $0.57$ & $-2.21$ & $-4.36$ & $-6.38$ & $-8.38$ &
   $2.22$ & $1.05$ & $-1.77$ & $-3.98$ & $-6.00$ & $-8.01$ 
\\
WJDTE 		& $0.66$ & $-1.06$ & $-3.04$ & $-5.03$ & $-7.03$ & $-9.03$ &
 $1.20$ & $-0.41$ & $-2.38$ & $-4.38$ & $-6.38$ & $-8.38$ &
  $1.48$ & $-0.06$ & $-2.01$ & $-4.01$ & $-6.01$ & $-8.01$ 
\\
JDTE-EP 		& $0.67$ & $-1.06$ & $-3.04$ & $-5.03$ & $-7.03$ & $-9.03$ &
 $1.20$ & $-0.41$ & $-2.38$ & $-4.38$ & $-6.38$ & $-8.38$ &
  $1.48$ & $-0.06$ & $-2.01$ & $-4.01$ & $-6.01$ & $-8.01$ 
\\
RJD			& $0.66$ & $-1.06$ & $-3.04$ & $-5.03$ & $-7.03$ & $-9.03$ &
 $1.20$ & $-0.41$ & $-2.38$ & $-4.38$ & $-6.38$ & $-8.38$ &
  $1.48$ & $-0.06$ & $-2.01$ & $-4.01$ & $-6.01$ & $-8.01$ 
\\
JD-MCG 		& $0.48$ & $-1.13$ & $-3.08$ & $-5.08$ & $-7.08$ & $-9.08$ &
 $0.99$ & $-0.62$ & $-2.41$ & $-4.40$ & $-6.40$ & $-8.40$ &
  $1.23$ & $-0.32$ & $-2.08$ & $-4.08$ & $-6.08$ & $-8.08$ 
\\
JD-MQN 		& $0.48$ & $-1.13$ & $-3.08$ & $-5.08$ & $-7.08$ & $-9.08$ &
 $0.99$ & $-0.62$ & $-2.41$ & $-4.40$ & $-6.40$ & $-8.40$ &
  $1.23$ & $-0.32$ & $-2.08$ & $-4.08$ & $-6.08$ & $-8.08$ 
 \\
\hline
 \end{tabular}
}
\caption{\label{fig:perf:SNR}{\bf Performance vs. noise level.\quad}
Median objective function values (top) and eigenvalue estimation error (bottom).}
\end{table*}

\subsection{Comparisons with state-of-the-art JD algorithms}

We compare our proposed methods with several recent state-of-the-art baselines from the literature: the algorithms WJDTE \cite{andre2020joint} and JDTE-EP \cite{cao2022joint} rely on simplifications of \eqref{functional} using Taylor expansions. 
The algorithm RJD from \cite{he2024randomized} considers a probabilistic approach which builds on successively diagonalizing random linear combinations of the set $\mathcal{A}$; we set here the number of combinations to $1000$.
Note that WJDTE outperformed many preceding approaches, see \cite{andre2020joint} for extensive comparisons; JDTE-EP is a recently proposed improvement of the method.
We also include our initial guess \eqref{equ:Uinit} as a baseline and denote it as $\mbox{eig}_{\sum\! A}$.

\noindent{\bf Performance versus iterations.\quad}
Fig.~\ref{fig:perf:iterations} plots the median objective function values as a function of iteration number for two different choices of $n$ and $K$ and SNR=$30$dB (RJD and $\mbox{eig}_{\sum\! A}$ are not shown because they are not iterative). The results lead to conclude that WJDTE and JDTE-EP stop their descent quickest, after less than 10 iterations, yet remain at an objective function value that is larger than that achieved by JD-MCG and JD-MQN: JD-MCG (JD-MQN) are at a lower objective function value than WJDTE and JDTE-EP latest at iteration 4 (iteration 3, respectively) and have converged to the same minimum with smaller objective function value of WJDTE and JDTE-EP shortly after.

\noindent{\bf Performance versus SNR.\quad}
Tab.~\ref{fig:perf:SNR} (top) shows the median of objective function values at the final iteration. As expected, all methods lead to better results for high SNR values. We observe that our proposed JD-MCG and JD-MQN yield the same and lowest objective function values - hence, find better minima of \eqref{functional2} than all baselines - consistently for all SNR values, and significantly so for small SNR. Second best are WJDTE and JDTE-EP, which show identical performance, followed by RJD and the initial guess $\mbox{eig}_{\sum\! A}$. The corresponding median estimation error for the eigenvalues $\Delta_k$ are shown in Tab.~\ref{fig:perf:SNR} (bottom) and lead to similar conclusions, except for very high SNR ($50$dB and $60$dB), in which all methods yield nearly identical solutions.

\begin{figure}[tb]
\centering
\includegraphics[width=0.45\linewidth]{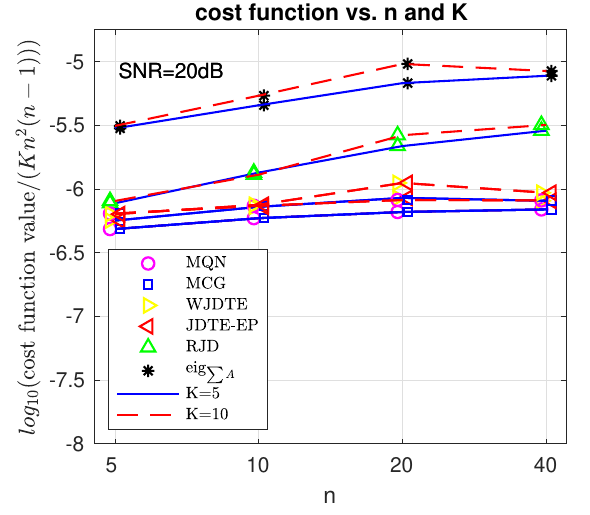}%
\includegraphics[width=0.45\linewidth]{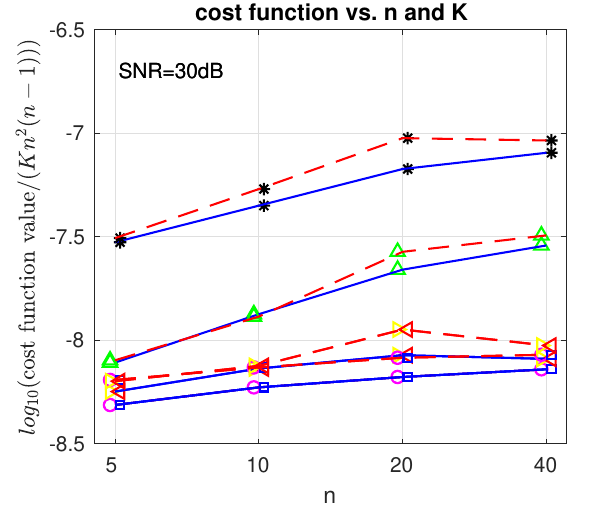}
\caption{\label{fig:perf:nk}{\bf\boldmath Performance comparisons vs. $n$ and $K$.\quad}
Median cost function values divided by $Kn^2(n-1)$ for various values for $n$ and $K$ and for   SNR=$20$dB (left) and SNR=$30$dB (right).}
\end{figure}

\noindent{\bf Performance versus $n$ and $K$.\quad}
A complementary performance comparison is provided in Fig.~\ref{fig:perf:nk}, where median cost function values (normalized by $Kn^2(n-1)$ for better visibility) are plotted for various values for $n$ and $K$. The results corroborate the conclusions from the previous paragraph, indicating that JD-MCG and JD-MQN reach better stationary points than all other baseline methods. 

\begin{figure}[tb]
\centering
\includegraphics[width=0.45\linewidth]{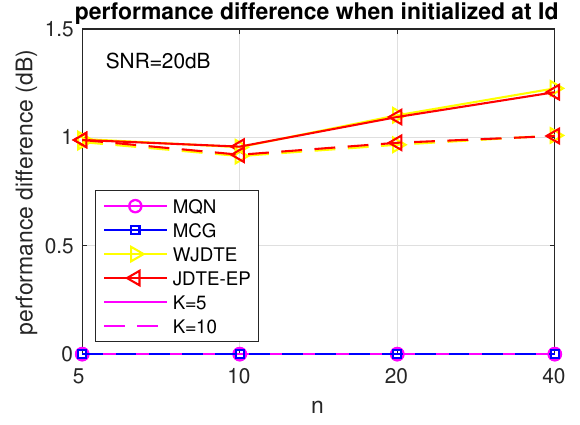}%
\includegraphics[width=0.45\linewidth]{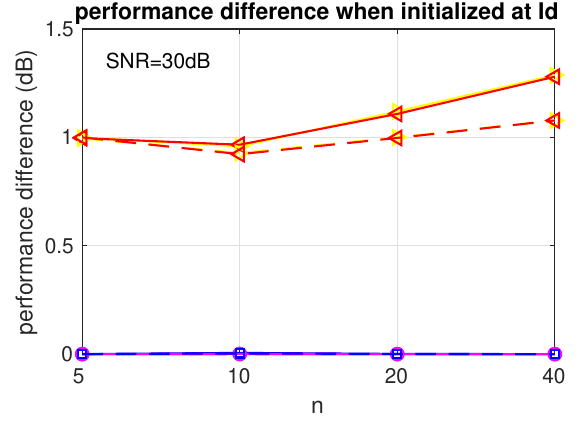}

\caption{\label{fig:perf:robustness}{\bf Robustness to initialization.\quad}
Ratio (in dB) of median cost function values reached when initializing at $U_{init}=I$ vs. \eqref{equ:Uinit} for various values for $n$ and $K$ and for   SNR=$30$dB (left) and SNR=$40$dB (right).}
\end{figure}

\noindent{\bf Robustness to initialization $U_{init}$.\quad}
To assess the robustness of our proposed algorithms and those of the iterative baseline algorithms WJDTE and JDTE-EP wrt.~initialization, we compare the objective function values they reach for initializations \eqref{equ:Uinit} (close to the minimum) and $U_{init}=I$ (far from the minimum).
Fig.~\ref{fig:perf:robustness} plots the objective function ratio for the different algorithms for several values for $K$, $n$ and SNR.
It can be observed that WJDTE and JDTE-EP overall yield solutions with $\approx1$dB larger objective function value when they are initialized at the identity matrix instead of  \eqref{equ:Uinit}.
On the contrary, our proposed JD-MCG and JD-MQN clearly converge to the same minima regardless of the initialization, indicating that our methods are more robust than WJDTE and JDTE-EP.

Overall, these results lead to conclude that our proposed algorithms are operational, effective and competitive. In comparison to the state of the art, they find better solutions in terms of objective function values at convergence, in particular for 
low SNR, and they are significantly more robust to initialization than the state of the art.

\section{Application to 3D harmonic retrieval}
\label{sec:3Dharmonic}
\begin{figure}[tb]
\centering
{\scriptsize\bf Data visualization for 3D harmonic retrieval}\\
\includegraphics[width=0.3\linewidth]{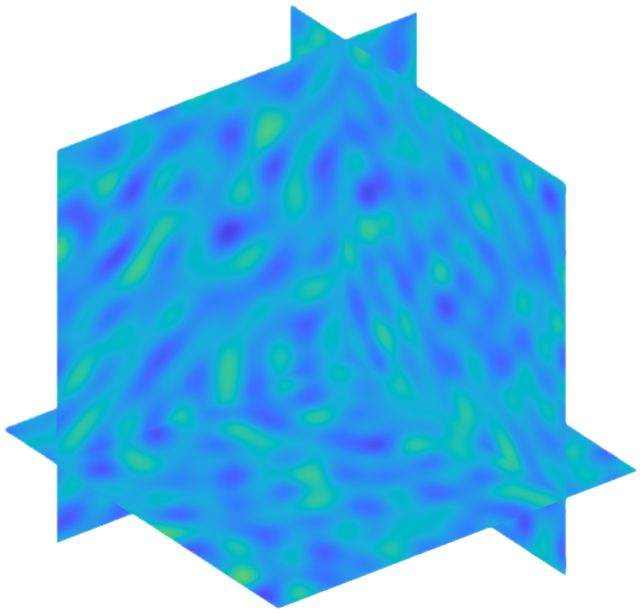}\\
\vskip-2mm{\small(a)}\\
\includegraphics[width=0.35\linewidth]{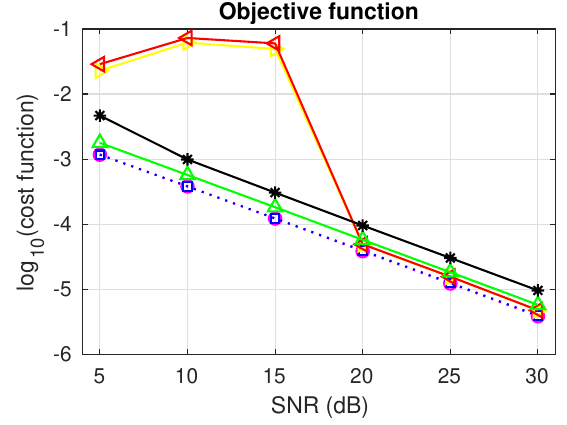}%
\includegraphics[width=0.35\linewidth]{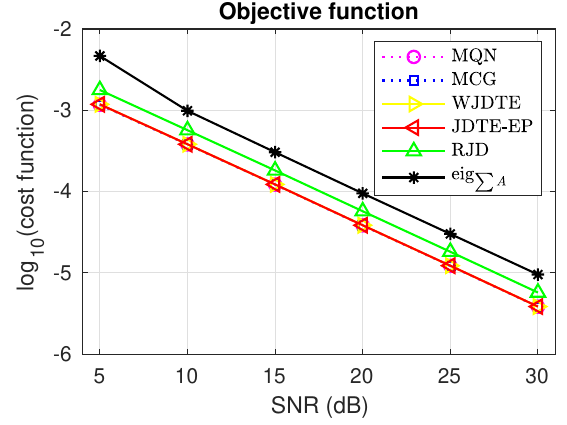}\\
\vskip-3mm{\small(b)}\\
\includegraphics[width=0.35\linewidth]{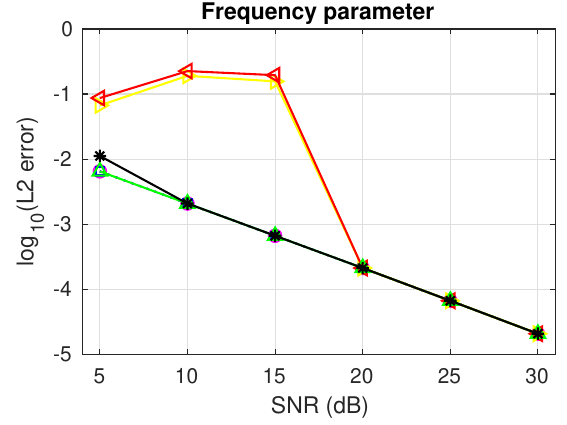}%
\includegraphics[width=0.35\linewidth]{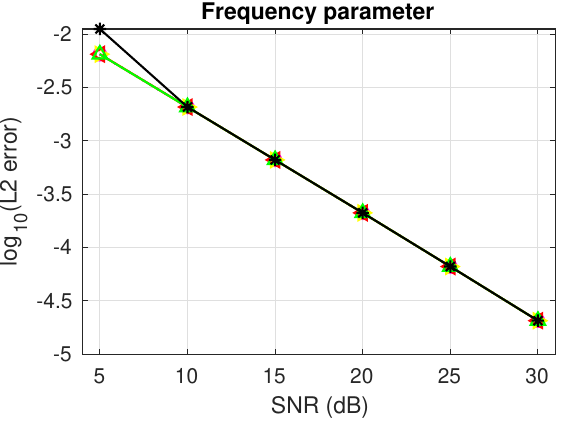}
\vskip-3mm{\small(c)}\\
\caption{\label{fig:perf:harmonicdata}{\bf Frequency estimation performance.}
Visualization of data cube for 3D harmonic retrieval (a);
median objective function value vs. SNR  (b) and median frequency estimation error (c) with initialization $U_{init}=I$ (left) and \eqref{equ:Uinit} (right).}
\end{figure}

We finally illustrate our algorithms for a 3D harmonic retrieval problem. The data consists of the sum of $K=27$ complex exponentials with random exponents and amplitudes and is sampled on a regular $17\times17\times17$ grid (see Fig.~\ref{fig:perf:harmonicdata} (a) for a visualization). The goal is to retrieve the harmonics from this data cube in circular Gaussian noise for various SNR. This is attempted using multidimensional ESPRIT \cite{sahnoun2017multidimensional,andersson2018esprit}, which leads to a JD problem of size $n=27$ and $K=3$.
Results for various levels of SNR are plotted in Fig.~\ref{fig:perf:harmonicdata}, in terms of median objective function value (b) and median frequency parameter estimation $\ell_2$ error (c), when initialized at $U_{init}=I$ (left column) and \eqref{equ:Uinit} (right column). The conclusions are similar to those obtained in the previous section: 
Our proposed JD-MCG and JD-MQN consistently yield lowest objective function values, regardless of the type of initialization, and for all SNR values. Nevertheless, in this problem, this only partially translates to better frequency estimates, and all methods yield similar results for the frequency parameters except at very low SNR values.
Note that WJDTE and JDTE-EP crucially rely on the use of good initializations to yield good performance in this application, otherwise they are lead astray, as the results for SNR levels below $20$dB for $U_{init}=I$ illustrate, for which both objective function values and frequency estimations are worse than the initial guess \eqref{equ:Uinit} denoted $\mbox{eig}_{\sum\! A}$.


\section{Conclusions}
\label{sec:conclusions}

In this work, we introduced novel first- and second-order descent algorithms for joint diagonalization (JD) of matrix collections, a fundamental task in fields such as signal processing, data analysis, and quantum information. Our approach overcomes limitations in existing JD methods by leveraging gradient and Hessian information in a computationally efficient manner. Specifically, we demonstrated how to evaluate the Hessian as a bilinear form or linear operator, bypassing the need to compute and store full second-order derivatives, which is critical for scalability in large-scale applications. Additionally, we incorporated and theoretically studied a multiplicative change of basis at each iteration, a modification that improves convergence speed and robustness.
The algorithms we proposed include gradient descent, conjugate gradient, and Quasi-Newton methods, each enhanced with a step-size selection strategy based on the Hessian’s structure. Our methods showed superior performance in numerical experiments, outperforming state-of-the-art techniques in solution quality, particularly in scenarios with high noise levels or large matrix dimensions.
Future work will focus on a deeper theoretical analysis of the properties of the objective function and further refinement of the algorithms for specific applications, such as high-dimensional harmonic retrieval in multidimensional signal processing.
Finally, to facilitate adoption, codes for our algorithms will be made public at the time of publication.

\bibliographystyle{IEEEtran}
\bibliography{JEVD}

\begin{thebibliography}{10}
\providecommand{\url}[1]{#1}
\csname url@samestyle\endcsname
\providecommand{\newblock}{\relax}
\providecommand{\bibinfo}[2]{#2}
\providecommand{\BIBentrySTDinterwordspacing}{\spaceskip=0pt\relax}
\providecommand{\BIBentryALTinterwordstretchfactor}{4}
\providecommand{\BIBentryALTinterwordspacing}{\spaceskip=\fontdimen2\font plus
\BIBentryALTinterwordstretchfactor\fontdimen3\font minus
  \fontdimen4\font\relax}
\providecommand{\BIBforeignlanguage}[2]{{%
\expandafter\ifx\csname l@#1\endcsname\relax
\typeout{** WARNING: IEEEtran.bst: No hyphenation pattern has been}%
\typeout{** loaded for the language `#1'. Using the pattern for}%
\typeout{** the default language instead.}%
\else
\language=\csname l@#1\endcsname
\fi
#2}}
\providecommand{\BIBdecl}{\relax}
\BIBdecl

\bibitem{jutten1991blind}
C.~Jutten and J.~H\'erault, ``Blind separation of sources, part i: An adaptive
  algorithm based on neuromimetic architecture,'' \emph{Signal Process.},
  vol.~24, no.~1, pp. 1--10, 1991.

\bibitem{cardoso1996equivariant}
J.-F. Cardoso and B.~Hvam~Laheld, ``Equivariant adaptive source separation,''
  \emph{IEEE T. Signal Process.}, vol.~44, no.~12, pp. 3017--3030, 1996.

\bibitem{pham2001joint}
D.~T. Pham, ``Joint approximate diagonalization of positive definite hermitian
  matrices,'' \emph{SIAM J. Matrix Analysis and Applications}, vol.~22, no.~4,
  pp. 1136--1152, 2001.

\bibitem{comon2010handbook}
P.~Comon and C.~Jutten, \emph{Handbook of Blind Source Separation: Independent
  Component Analysis and Applications}.\hskip 1em plus 0.5em minus 0.4em\relax
  Elsevier Science, 2010.

\bibitem{cardoso1993blind}
J.-F. Cardoso and A.~Souloumiac, ``Blind beamforming for non-gaussian
  signals,'' in \emph{IEE Proc. F (Radar and Signal Process.)}, vol. 140,
  no.~6, 1993, pp. 362--370.

\bibitem{comon1994independent}
P.~Comon, ``Independent component analysis, a new concept?'' \emph{Signal
  Process.}, vol.~36, no.~3, pp. 287--314, 1994.

\bibitem{ziehe2004fast}
A.~Ziehe, P.~Laskov, G.~Nolte, and K.-R. M{\"u}ller, ``A fast algorithm for
  joint diagonalization with non-orthogonal transformations and its application
  to blind source separation,'' \emph{J. Machine Learning Research}, vol.~5,
  no. Jul, pp. 777--800, 2004.

\bibitem{flury1984common}
B.~N. Flury, ``Common principal components in k groups,'' \emph{J. American
  Statistical Association}, vol.~79, no. 388, pp. 892 -- 898, 1984.

\bibitem{roemer2013semi}
F.~Roemer and M.~Haardt, ``A semi-algebraic framework for approximate cp
  decompositions via simultaneous matrix diagonalizations (secsi),''
  \emph{Signal Process.}, vol.~93, no.~9, pp. 2722--2738, 2013.

\bibitem{luciani2014canonical}
X.~Luciani and L.~Albera, ``Canonical polyadic decomposition based on joint
  eigenvalue decomposition,'' \emph{Chemometrics and Intelligent Laboratory
  Systems}, vol. 132, pp. 152--167, 2014.

\bibitem{xu2024testing}
Y.~Xu, M.-C. D{\"u}ker, and D.~S. Matteson, ``Testing simultaneous
  diagonalizability,'' \emph{J. American Statistical Association}, vol. 119,
  no. 546, pp. 1513--1525, 2024.

\bibitem{haardt2008higher}
M.~Haardt, F.~Roemer, and G.~Del~Galdo, ``Higher-order svd-based subspace
  estimation to improve the parameter estimation accuracy in multidimensional
  harmonic retrieval problems,'' \emph{IEEE T. Signal Process.}, vol.~56,
  no.~7, pp. 3198--3213, 2008.

\bibitem{andersson2015method}
F.~Andersson, M.~Carlsson, J.-Y. Tourneret, and H.~Wendt, ``A method for 3d
  direction of arrival estimation for general arrays using multiple
  frequencies,'' in \emph{IEEE Int. Workshop on Computational Advances in
  Multi-Sensor Adaptive Process. (CAMSAP)}.\hskip 1em plus 0.5em minus
  0.4em\relax IEEE, 2015, pp. 325--328.

\bibitem{sahnoun2017multidimensional}
S.~Sahnoun, K.~Usevich, and P.~Comon, ``Multidimensional esprit for damped and
  undamped signals: Algorithm, computations, and perturbation analysis,''
  \emph{IEEE T. Signal Process.}, vol.~65, no.~22, pp. 5897--5910, 2017.

\bibitem{andersson2018esprit}
F.~Andersson and M.~Carlsson, ``Esprit for multidimensional general grids,''
  \emph{SIAM J. Matrix Analysis and Applications}, vol.~39, no.~3, pp.
  1470--1488, 2018.

\bibitem{behmanesh2021geometric}
M.~Behmanesh, P.~Adibi, J.~Chanussot, C.~Jutten, and S.~M.~S. Ehsani,
  ``Geometric multimodal learning based on local signal expansion for joint
  diagonalization,'' \emph{IEEE T. Signal Process.}, vol.~69, pp. 1271--1286,
  2021.

\bibitem{burel2022joint}
G.~Burel, H.~Pillin, P.~Baird, E.-H. Baghious, and R.~Gautier, ``Joint
  eigenvalue decomposition for quantum information theory and processing,''
  \emph{Matrix Theory-Classics and Advances}, 2022.

\bibitem{pham2001blind}
D.-T. Pham and J.-F. Cardoso, ``Blind separation of instantaneous mixtures of
  nonstationary sources,'' \emph{IEEE T. Signal Process.}, vol.~49, no.~9, pp.
  1837--1848, 2001.

\bibitem{yeredor2002non}
A.~Yeredor, ``Non-orthogonal joint diagonalization in the least-squares sense
  with application in blind source separation,'' \emph{IEEE T. Signal
  Process.}, vol.~50, no.~7, pp. 1545--1553, 2002.

\bibitem{tichavsky2008fast}
P.~Tichavsky and A.~Yeredor, ``Fast approximate joint diagonalization
  incorporating weight matrices,'' \emph{IEEE T. Signal Process.}, vol.~57,
  no.~3, pp. 878--891, 2008.

\bibitem{souloumiac2009nonorthogonal}
A.~Souloumiac, ``Nonorthogonal joint diagonalization by combining givens and
  hyperbolic rotations,'' \emph{IEEE T. Signal Process.}, vol.~57, no.~6, pp.
  2222--2231, 2009.

\bibitem{ablin2018faster}
P.~Ablin, J.-F. Cardoso, and A.~Gramfort, ``Faster independent component
  analysis by preconditioning with hessian approximations,'' \emph{IEEE T.
  Signal Process.}, vol.~66, no.~15, pp. 4040--4049, 2018.

\bibitem{azzouz2023generalized}
M.~Azzouz, A.~Mesloub, K.~Abed-Meraim, and A.~Belouchrani, ``Generalized
  unitary joint diagonalization algorithm based on approximate givens
  rotations,'' \emph{IEEE Signal Process. Letters}, 2023.

\bibitem{luciani2015joint}
X.~Luciani and L.~Albera, ``Joint eigenvalue decomposition of non-defective
  matrices based on the lu factorization with application to ica,'' \emph{IEEE
  T.. Signal Process.}, vol.~63, no.~17, pp. 4594--4608, 2015.

\bibitem{mesloub2018efficient}
A.~Mesloub, A.~Belouchrani, and K.~Abed-Meraim, ``Efficient and stable joint
  eigenvalue decomposition based on generalized givens rotations,'' in
  \emph{Proc. European Signal Process. Conf. (EUSIPCO)}, Rome, Italy, 2018.

\bibitem{iferroudjene2009new}
R.~Iferroudjene, K.~A. Meraim, and A.~Belouchrani, ``A new jacobi-like method
  for joint diagonalization of arbitrary non-defective matrices,''
  \emph{Applied Math. and Comput.}, vol. 211, no.~2, pp. 363--373, 2009.

\bibitem{gong2012complex}
X.-F. Gong, K.~Wang, and Q.-H. Lin, ``Complex non-orthogonal joint
  diagonalization with successive givens and hyperbolic rotations,'' in
  \emph{Proc. Int. Conf. Acoustics Speech and Signal Process. (ICASSP)}, Kyoto,
  Japan, 2012.

\bibitem{he2024randomized}
H.~He and D.~Kressner, ``Randomized joint diagonalization of symmetric
  matrices,'' \emph{SIAM J. Matrix Analysis and Applications}, vol.~45, no.~1,
  pp. 661--684, 2024.

\bibitem{andre2020joint}
R.~Andr{\'e}, X.~Luciani, and E.~Moreau, ``Joint eigenvalue decomposition
  algorithms based on first-order {T}aylor expansion,'' \emph{IEEE T. Signal
  Process.}, vol.~68, pp. 1716--1727, 2020.

\bibitem{cao2022joint}
Q.-C. Cao, G.-H. Cheng, and E.~Moreau, ``The joint eigenvalue decomposition
  algorithm based on first-order taylor expansion via the exterior penalty
  function method,'' \emph{Signal Process.}, vol. 200, p. 108644, 2022.

\bibitem{hori1999joint}
G.~Hori, ``Joint diagonalization and matrix differential equations,'' in
  \emph{Proc. Int. Symp. Nonlinear Theory and Applications (NOLTA)}, Hawaii,
  USA, 1999.

\bibitem{troedsson2024joint}
E.~Troedsson, D.~Falkowski, C.-F. Lidgren, H.~Wendt, and M.~Carlsson, ``On
  joint eigen-decomposition of matrices,'' \emph{arXiv preprint
  arXiv:2409.10292}, 2024.

\bibitem{ErikEUSIPCO2024}
E.~Troedsson, M.~Carlsson, and H.~Wendt, ``On gradient based descent algorithms
  for joint diagonalization of matrices,'' in \emph{Proc. European Signal
  Process. Conf. (EUSIPCO)}, Lyon, France, Aug. 2024.

\bibitem{Yeredor2005}
A.~Yeredor, ``On using exact joint diagonalization for noniterative approximate
  joint diagonalization,'' \emph{IEEE Signal Process. Letters}, vol.~12, no.~9,
  pp. 645--648, 2005.

\bibitem{daniel1967conjugate}
J.~W. Daniel, ``The conjugate gradient method for linear and nonlinear operator
  equations,'' \emph{SIAM J. Numerical Analysis}, vol.~4, no.~1, pp. 10--26,
  1967.

\bibitem{dai1999nonlinear}
Y.~H. Dai and Y.~Yuan, ``A nonlinear conjugate gradient method with a strong
  global convergence property,'' \emph{SIAM J. Optimization}, vol.~10, no.~1,
  pp. 177--182, 1999.

\bibitem{fletcher1964function}
R.~Fletcher and C.~M. Reeves, ``Function minimization by conjugate gradients,''
  \emph{The computer journal}, vol.~7, no.~2, pp. 149--154, 1964.

\bibitem{hestenes1952methods}
M.~R. Hestenes and E.~Stiefel, ``Methods of conjugate gradients for solving
  linear systems,'' \emph{J. Research of the National Bureau of Standards},
  vol.~49, pp. 409--435, 1952.

\bibitem{polak1969note}
E.~Polak and G.~Ribiere, ``\BIBforeignlanguage{fr}{Note sur la convergence de
  m\'ethodes de directions conjugu\'ees},'' \emph{\BIBforeignlanguage{fr}{Revue
  fran\c{c}aise d'informatique et de recherche op\'erationnelle. S\'erie
  rouge}}, vol.~3, no.~R1, pp. 35--43, 1969.

\bibitem{william2006survey}
W.~W. Hager and H.~Zhang, ``A survey of nonlinear conjugate gradient methods,''
  \emph{Pacific Journal of Optimization}, vol.~2, no.~1, pp. 35--58, 2006.

\bibitem{gilbert1992global}
J.~C. Gilbert and J.~Nocedal, ``Global convergence properties of conjugate
  gradient methods for optimization,'' \emph{SIAM J. Optimization}, vol.~2,
  no.~1, pp. 21--42, 1992.

\end{thebibliography}

 \end{document}